\documentclass[moor]{informs1}

\usepackage{natbib}
 \NatBibNumeric
 \def\bibfont{\small}%
 \bibpunct[, ]{[}{]}{,}{n}{}{,}%

\usepackage[colorlinks=true,breaklinks=true,bookmarks=true,urlcolor=blue,
     citecolor=blue,linkcolor=blue,bookmarksopen=false,draft=false]{hyperref}

\def\EMAIL#1{\href{mailto:#1}{#1}}
\def\URL#1{\href{#1}{#1}}

\TheoremsNumberedThrough
\EquationsNumberedThrough
\MANUSCRIPTNO{N/A} 

\newcommand{\proofbeg}{\begin{proof}{Proof.}}
\newcommand{\proofof}[1]{\begin{proof}{#1}}
\newcommand{\proofend}{\hfill $\square$ \end{proof} \smallskip}

\newcommand{\citemor}[1]{\citeauthor{#1}~\cite{#1}}
\newcommand{\citemorpar}[1]{\citeauthor{#1}~\cite{#1}}
\newcommand{\citemoradd}[2]{\citeauthor{#1}~\cite[#2]{#1}}

\newenvironment{rema}{\smallskip\begin{remark}}{\end{remark}\smallskip}
\newenvironment{defi}{\smallskip\begin{definition}}{\end{definition}\smallskip}

\newenvironment{MORproof}{\smallskip\proofbeg}{\proofend}

\usepackage{bigints}
\usepackage[ruled,norelsize]{algorithm2e}
\usepackage{undertildenu} 

\newcommand{\R}{\mathbb{R}}
\newcommand{\ind}{\mathbb{I}}
\newcommand{\E}{\mathbb{E}}
\newcommand{\KL}{\mathrm{KL}}
\newcommand{\Kinf}{\mathcal{K}_{\inf}}
\newcommand{\kl}{\mathrm{kl}}
\renewcommand{\P}{\mathbb{P}}
\renewcommand{\geq}{\geqslant}
\renewcommand{\leq}{\leqslant}
\newcommand{\m}{\mathrm{m}}
\newcommand{\mset}{\mathcal{M}\bigl([0,M]\bigr)}
\newcommand{\msetfin}{\mathcal{M}_{\mathrm{fin}}\bigl([0,M]\bigr)}
\newcommand{\eps}{\varepsilon}
\renewcommand{\d}{\mathrm{d}}
\newcommand{\algo}{\texttt}
\newcommand{\cA}{\mathcal{A}}
\newcommand{\cK}{\mathcal{K}}
\newcommand{\cD}{\mathcal{D}}
\newcommand{\cB}{\mathcal{B}}
\newcommand{\cG}{\mathcal{G}}
\newcommand{\uDelta}{\underline{\Delta}}

\renewcommand{\O}{\mathrm{O}}
\newcommand{\nuun}{\underline{\underline{\nu}}}
\newcommand{\nuzero}{\utilde{\nu}}
\newcommand{\cW}{\mathcal{W}}
\newcommand{\tw}{\widetilde{w}}
\newcommand{\cN}{\mathcal{N}}
\newcommand{\unu}{\underline{\nu}}

\begin{document}

\RUNAUTHOR{Garivier, M{\'e}nard, Stoltz}
\RUNTITLE{The True Shape of Regret in Bandit problems}
\TITLE{Explore First, Exploit Next: \\ The True Shape of Regret in Bandit Problems}

\ARTICLEAUTHORS{%
\AUTHOR{Aur{\'e}lien Garivier}
\AFF{IMT: Universit{\'e} Paul Sabatier -- CNRS, Toulouse, France, \EMAIL{aurelien.garivier@math.univ-toulouse.fr}, \URL{http://www.math.univ-toulouse.fr/~agarivie/}}
\AUTHOR{Pierre M{\'e}nard}
\AFF{IMT: Universit{\'e} Paul Sabatier -- CNRS, Toulouse, France}
\AUTHOR{Gilles Stoltz}
\AFF{GREGHEC: HEC Paris -- CNRS, Jouy-en-Josas, France, \EMAIL{stoltz@hec.fr},
\URL{http://stoltz.perso.math.cnrs.fr}}
}

\ABSTRACT{%
We revisit lower bounds on the regret in the case of multi-armed bandit problems.
We obtain non-asymptotic, distribution-dependent bounds and provide simple proofs based only on well-known properties of Kullback-Leibler divergences.
These bounds show in particular that in the initial phase the regret grows almost linearly, and that the well-known logarithmic growth of the regret only holds in a final phase.
The proof techniques come to the essence of the information-theoretic arguments used and they involve no unnecessary complications.
}%

\MSCCLASS{Primary: 68T05, 62L10; Secondary: 62B10, 94A17, 94A20}
\ORMSCLASS{Primary: computer science: artificial intelligence (learning and adaptive systems);
statistics: sequential methods (sequential analysis); secondary:
statistics: sufficiency and information (information-theoretic topics);
information and communication, circuits: communication, information (measures of information,
entropy; sampling theory)}
\HISTORY{Submitted June 10, 2016; revised March 3, 2017; accepted December 8, 2017}

\maketitle

\section{Introduction.}

\begin{figure}
\center
\begin{tabular}{lr}
\includegraphics[width=0.47\textwidth]{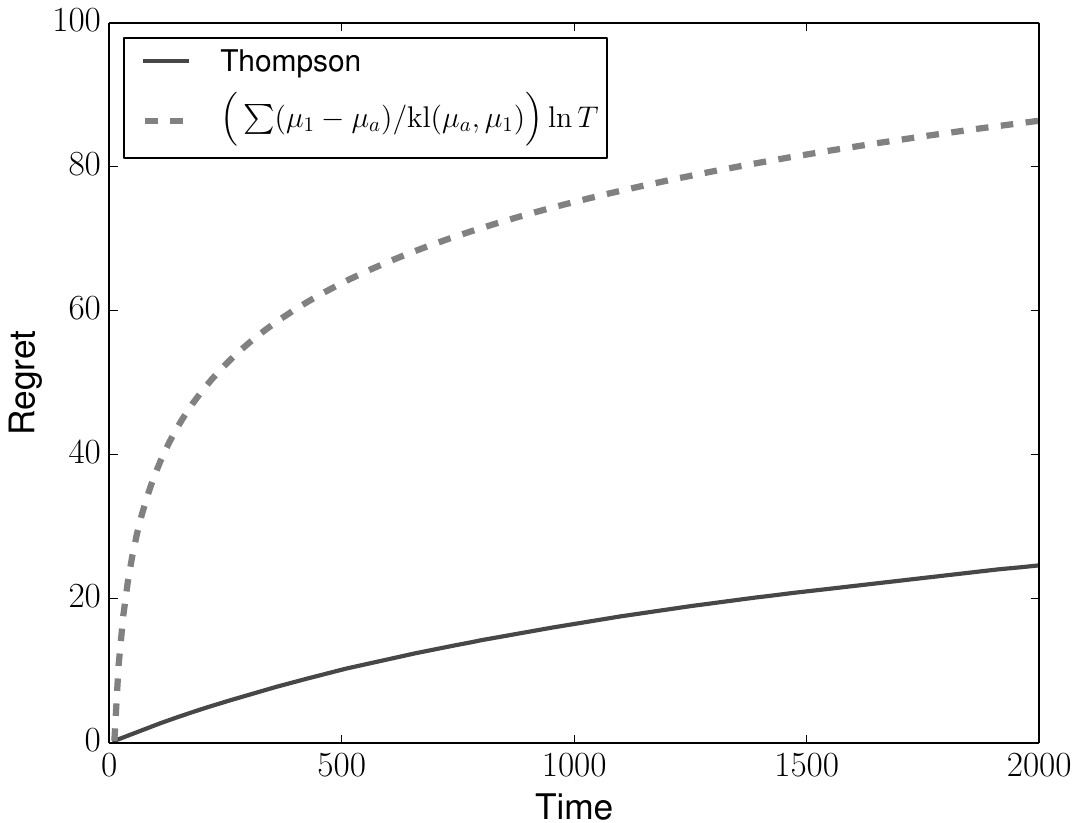} &
\includegraphics[width=0.47\textwidth]{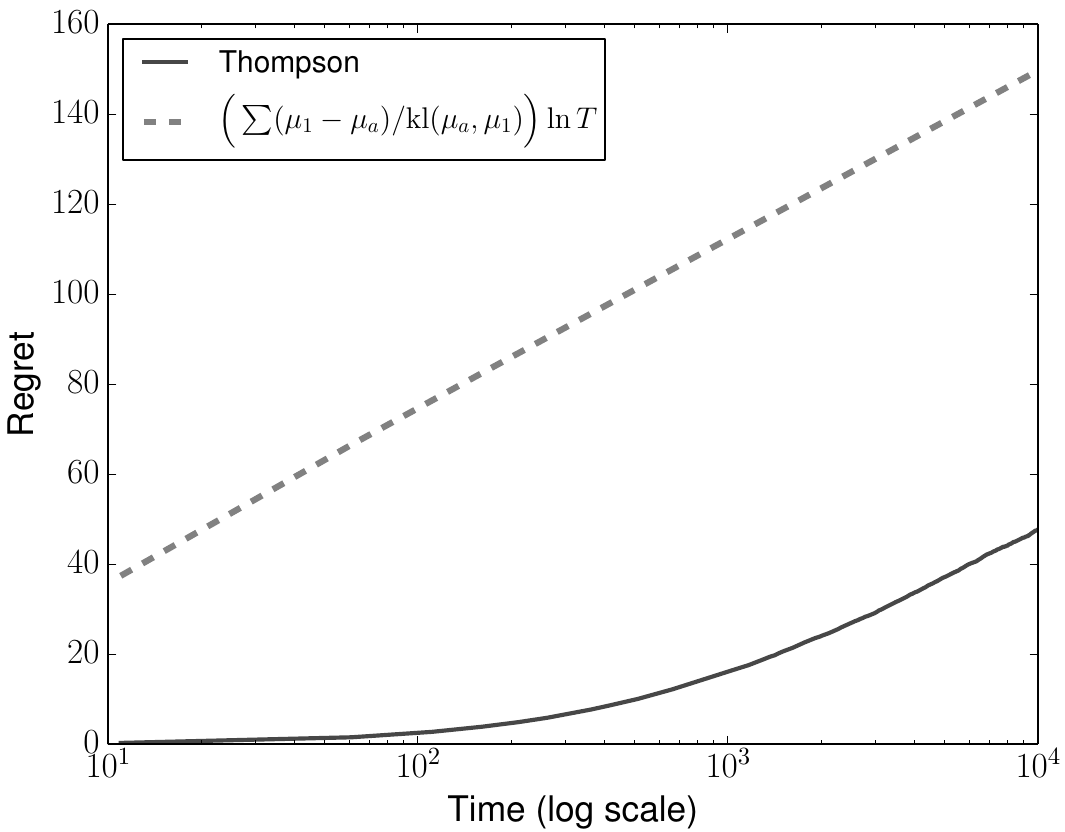}
\end{tabular}
\caption{\label{fig:draws}
Expected regret of \citemor{Th33} Sampling (\emph{blue, solid} line) on a Bernoulli bandit problem with parameters $(\mu_a)_{1 \leq a \leq 6} =(0.05, \, 0.04, \, 0.02, \, 0.015, \, 0.01, \, 0.005)$; expectations are approximated over $500$ runs. \smallskip \newline
Versus the \citemor{LaRo85} lower bound (\emph{red, dotted} line) for a Bernoulli model;
here $\kl$ denotes the Kullback-Leibler divergence~\eqref{eq:kl} between Bernoulli distributions.
\smallskip \newline
\emph{Left}: the shape of regret is not logarithmic at first, rather linear. \newline
\emph{Right}: the asymptotic lower bound is out of reach unless $T$ is extremely large.}
\end{figure}

After the works of~\citemor{LaRo85} and \citemor{BuKa96}, it is widely admitted that the growth of the cumulative regret in a bandit problem is a logarithmic function of time, multiplied by a sum of terms involving Kullback-Leibler divergences.
The asymptotic nature of the lower bounds, however, appears clearly in numerical experiments, where the logarithmic shape is not to be observed on small horizons (see Figure~\ref{fig:draws}, left). Even on larger horizons, the second-order terms keep a large importance, which causes the regret of some algorithms to remain way \emph{below} the ``lower bound'' on any experimentally visible horizon (see Figure~\ref{fig:draws}, right; see also~\citemor{ETC}).

\paragraph{First contribution: a folk result made rigorous.}
It seems to be a folk result (or at least, a widely believed result) that the regret should be linear in an initial phase of a bandit problem.
However, all references that we were pointed out exhibit such a linear behavior only for limited
bandit settings; we discuss them below, in the section about literature review. We are the first
to provide linear distribution-dependent lower bounds for small horizons that hold for general
bandit problems, with no restriction on the shape or on the expectations of the distributions over the arms.

Thus we may draw a more precise picture of the behavior of the regret in any bandit problem.
Indeed, our bounds show the existence of three successive phases: an initial linear phase, when all the arms are essentially drawn uniformly; a transition phase, when the number of observations becomes sufficient to perceive differences; and the final phase, when the distributions associated with all the arms are known with high confidence and when the new draws are just confirming the identity of the best arms with higher and higher degree of confidence (this is the famous logarithmic phase).
This last phase may often be out of reach in applications, especially when the number of arms is large.

\paragraph{Second contribution: a generic tool for proving distribution-dependent bandit lower bounds.}
On the technical side, we provide simple proofs, based on the fundamental
information-theoretic inequality~\eqref{eq:funda} stated in Section~\ref{sec:fundamental},
which generalizes and simplifies previous approaches
based on explicit changes of measures. In particular, we are able to re-derive the
asymptotic distribution-dependent lower bounds of~\citemor{LaRo85}, \citemor{BuKa96}
and \citemor{CoKa} in a few lines.
This may perhaps be one of the most striking contributions of this paper.
As a final set of results, we offer non-asymptotic versions of these
lower bounds for large horizons, and exhibit the optimal order of magnitude
of the second-order term in the regret bound, namely, $-\ln(\ln T)$.

The proof techniques come to the essence of the arguments used so far in the literature
and they involve no unnecessary complications; they only rely on well-known properties of
Kullback-Leibler divergences.

\subsection{Setting.}
\label{sec:setting}

We consider the simplest case of a stochastic bandit problem,
with finitely many arms indexed by $a \in \{1,\ldots,K\}$.
Each of these arms is associated with an unknown probability distribution $\nu_a$ over~$\R$.
We assume that each $\nu_a$ has a well-defined expectation
and call $\unu = (\nu_a)_{a = 1,\ldots,K}$ a bandit problem.

At each round $t \geq 1$, the player pulls the arm $A_t$ and gets a real-valued reward $Y_t$ drawn
independently at random
according to the distribution $\nu_{A_t}$.
This reward is the only piece of information
available to the player.

\paragraph{Strategies.}
A strategy $\psi$ associates an arm with the information gained in the past, possibly
based on some auxiliary randomization; without loss of generality, this auxiliary randomization
is provided by a sequence $U_0,U_1,U_2,\ldots$ of independent and identically distributed
random variables, with common distribution the uniform distribution over $[0,1]$.
Formally, a strategy is a sequence $\psi = (\psi_t)_{t \geq 0}$ of measurable functions, each of which
associates with the said past information, namely,
\[
I_t = \bigr( U_0,\,Y_1,U_1,\,\ldots,\,Y_t,U_t \bigr)\,,
\]
an arm $\psi_t(I_t) = A_{t+1} \in \{1,\ldots,K\}$, where $t \geq 0$. The initial information
reduces to $I_0 = U_0$ and the first arm is $A_1 = \psi_0(U_0)$. The auxiliary
randomization is conditionally independent of the sequence of rewards in the following sense:
for $t \geq 1$, the randomization $U_{t}$ used to pick $A_{t+1}$
is independent of $I_{t-1}$ and $Y_t$.

\paragraph{Regret.}
A typical measure of the performance of a strategy is given by its regret.
To recall its definition, we denote by $E(\nu_a) = \mu_a$ the expected payoff of arm $a$
and by $\Delta_a$ its gap to an optimal arm:
\[
\mu^\star = \max_{a=1,\ldots,K} \mu_a \qquad \mbox{and} \qquad
\Delta_a = \mu^\star - \mu_a\,.
\]
The number of times an arm $a$ is pulled until round $T$ by a strategy $\psi$ is referred to
as
\[
N_{\psi,a}(T) = \sum_{t=1}^T \ind_{\{ A_t = a \}} = \sum_{t=1}^T \ind_{\{ \psi_{t-1}(I_{t-1}) = a \}} \,.
\]
The expected regret of a strategy $\psi$ equals, by the tower rule (see details below),
\begin{equation}
\label{eq:towerrule}
R_{\psi,\unu,T} = T\mu^\star - \E_{\unu}\!\!\left[ \sum_{t=1}^T Y_t \right]
= \E_{\unu}\!\!\left[ \sum_{t=1}^T \bigl( \mu^\star - \mu_{A_t} \bigr) \right]
= \sum_{a=1}^K \Delta_a \, \E_{\unu} \! \bigl[ N_{\psi,a}(T) \bigr]\,.
\end{equation}
In the equation above, the notation $\E_{\unu}$ refers to the expectation
associated with the bandit problem $\unu = (\nu_a)_{a = 1,\ldots,K}$;
it is made formal in Section~\ref{sec:fundamental}.

To show~\eqref{eq:towerrule}, we use that by the definition of the bandit setting,
the distribution of the obtained payoff $Y_t$ only depends on the chosen arm
$A_t$ and is independent from the past random draws of the $Y_1,\ldots,Y_{t-1}$.
More precisely, conditionally on $A_t$, the distribution of $Y_t$ is $\nu_{A_t}$ so that
\[
\E_{\unu}\!\bigl[ Y_t \,|\, A_t \bigr] = \mu_{A_t}\,, \qquad \mbox{thus} \qquad
\E_{\unu}[Y_t] = \E_{\unu} \!\Bigl[ \E_{\unu}\!\bigl[ Y_t \,|\, A_t \bigr] \Bigr] = \E_{\unu} \! \bigl[ \mu_{A_t} \bigr]\,,
\]
where we used the tower rule for the second set of equalities.

\subsection{The general asymptotic lower bound: a quick literature review.}

We consider a bandit model $\mathcal{D}$, i.e., a collection of possible
distributions $\nu_a$ associated with the arms.
(That is, $\mathcal{D}$ is a subset of the set of all possible distributions
over $\R$ with an expectation.)
\citemor{LaRo85} and later \citemor{BuKa96} exhibited asymptotic lower
bounds and matching asymptotic upper bounds on the
normalized regret $R_{\psi,\unu,T}/\ln T$,
respectively in a one-parameter case and in a more general, multi-dimensional parameter case,
under mild conditions on $\mathcal{D}$.
We believe that
the extension of these bounds to any, even non-parametric, model was a known or at least conjectured result
(see, for instance, the introduction of~\citemorpar{klucb}). It turns out that
recently, \citemor{CoKa} provided a clear non-parametric statement, though under additional mild
conditions on the model $\cD$, which, as we will see, are not needed.

In the sequel,
we denote by $\KL$ the Kullback-Leibler divergence between
two probability distributions. We also recall that we denoted by $E$ the expectation operator
(that associates with each distribution its expectation).

To state the bound for the case of an arbitrary model $\mathcal{D}$, we will
use the following key quantity $\Kinf$ introduced by \citemoradd{BuKa96}{quantity (3)--(b) on page 125}.

\paragraph{The key quantity $\Kinf$.}
For any given $\nu_a \in \mathcal{D}$ and
any real number $x$,
\[
\Kinf(\nu_a,x,\cD) = \inf \Bigl\{ \KL(\nu_a,\nu'_a) : \ \ \nu'_a \in \mathcal{D} \ \ \mbox{and} \ \
E(\nu'_a) > x \Bigr\}\,;
\]
by convention, the infimum of the empty set equals $+\infty$.
When the considered strategy is uniformly fast convergent in the sense of Definition~\ref{eq:consistency} (stated later in this paper),
then, for any suboptimal arm $a$,
\begin{equation}
\label{eq:obj}
\liminf_{T \to \infty} \,\, \frac{\E_{\unu} \! \bigl[ N_{\psi,a}(T) \bigr]}{\ln T}
\geq \frac{1}{\Kinf(\nu_a,\mu^\star,\cD)}\,.
\end{equation}
Note that by the convention on the infimum of the empty set,
this lower bound is void as soon as there exists no $\nu'_a \in \mathcal{D}$
such that $E(\nu'_a) > \mu^\star$.

\paragraph{Previous partial simplifications of the proof of~\eqref{eq:obj}.}
We re-derive the above bound in a few lines in Section~\ref{sec:rederiv-distrfree}.

There had been recent attempts to clarify the exposition of the proof of this lower bound,
together with the desire of dropping the mild conditions that were still needed so far on the model $\cD$. We first
mention that \citemor{CoKa} provided a more general and streamlined approach
than the original expositions by \citemor{LaRo85} and \citemor{BuKa96}.

The case of Bernoulli models was discussed in~\citemor{Bu10} and \citemor{BuCB12}.
Only assumptions of uniform fast convergence of the strategies are required (see Definition~\ref{eq:consistency}) and the associated
proof follows the original proof technique, by performing first an explicit change of measure
and then applying some Markov--Chernoff bounding. More recently,
\citemoradd{CJ15}{Section~2.2} presented a proof (only in the Bernoulli case)
not relying on any explicit change of measure but with many additional technicalities
with respect to our exposition, including some Markov bounding of well-chosen events.
We have been referred to this PhD dissertation only recently,
after completing the present paper.

As far as general bandit models are concerned, we may cite \citemoradd{KaCaGa16}{Appendix~B}:
they deal with the case of any model $\mathcal{D}$
but with the restriction that only bandit problems $\unu$
with a unique optimal arm should be considered. They still use both an explicit change of measure --to prove
the chain-rule equality in~\eqref{eq:funda}-- and then apply as well
some Markov--Chernoff bounding to the probability of well-chosen events.
With a different aim, \citemor{CoPr14} presented similar arguments.

We also wish to mention the contribution of~\citemor{WuGySz15bi},
though their focus and aim are radically different.
With respect to some aspects, their setting and goal is wider or more general:
they developed non-asymptotic problem-dependent lower bounds on the regret of any algorithm, in the case
of more general limited feedback models than just the simplest case of multi-armed
bandit problems. Their lower bounds can recover the asymptotic bounds
of \citemor{BuKa96}, but only up to a constant factor as they acknowledge in their contribution.
These lower bounds are in terms of uniform upper bounds on the regret of the considered strategies,
which is in contrast with the lower bounds we develop in Section~\ref{sec:smallT}. Therein, we need
some assumptions on the strategies --extremely mild ones, though: some minimal symmetry-- and do not need
their regret to be bounded from above.
However, the main difference with respect to this reference is that
its focus is limited to specific bandit models, namely Gaussian bandits models, while
\citemor{BuKa96} and the present paper do not impose such a restriction on the bandit model.

\subsection{Other bandit lower bounds: a brief literature review.}
\label{sec:other-blb}
In this paper, we are mostly interested in general distribution-dependent lower bounds,
that hold for all bandit problems, just like~\eqref{eq:obj}. We do target generality.
This is in contrast with many earlier lower bounds in the multi-armed bandit setting,
which are rather of the following form, which we will refer to as (well-chosen):
\begin{quote}
``\emph{There exists some \emph{well-chosen}, difficult bandit problem such that all strategies suffer
a regret larger than} [...].'' \hfill (well-chosen)
\end{quote}
Specific examples and pointers for this kind of bounds are given below.
An interesting variation is provided by~\citemoradd{MT04}{Theorem~10},
who state that for all strategies, there exists some well-chosen, difficult Bernoulli bandit problem
such that the regret is linear at first and then, logarithmic.

On the contrary, we will issue statements of the following form, which we will refer to as (all):
\begin{quote}
``\emph{For \emph{all} bandit problems, all (reasonable) strategies suffer a regret larger than} (...).''
\hfill (all)
\end{quote}
Sometimes, but not always, we will have to impose some mild restrictions on the considered
strategies (like some minimal symmetry, or some notion of uniform fast convergence); this is what we
mean by requiring the strategies to be ``reasonable''.

We discuss briefly below two other sets of regret lower bounds.
We are pleased to mention that our fundamental inequality was already
used in at least one subsequent article, namely by~\citemor{ETC},
to prove in a few lines matching lower bounds for a refined
analysis of explore-then-commit strategies.

\paragraph{The distribution-free lower bound.}
This inequality states that for the model $\mathcal{D} = \mathcal{M}\bigl([0,1]\bigr)$ of all
probability distributions over $[0,1]$,
for all strategies $\psi$,
for all $T \geq 1$ and all $K \geq 2$,
\begin{equation}
\label{eq:distrfree}
\sup_{\unu} R_{\psi,\unu,T} \geq \frac{1}{20} \min\Bigl\{ \sqrt{KT}, \, T \Bigr\}\,;
\end{equation}
see~\citemor{AuCBFrSc02}, \citemor{CBLu06}, and for two-armed bandits, \citemor{KuLu00}.
We re-derive the above bound in Section~\ref{sec:rederiv-distrfree} of the appendix.
This re-derivation follows the very same proof scheme as in the original proof; the
only difference is that some steps (e.g., the use of chain-rule equality for Kullback-Leibler
divergences) are implemented separately as parts of the proof of our general
inequality~\eqref{eq:funda}.
In particular, the well-chosen difficult bandit problems used to prove this bound
are composed of Bernoulli distributions with parameters $1/2$ and $1/2+\varepsilon$,
where $\varepsilon$ is carefully tuned according to the values of $T$ and $K$.
This bound therefore rather falls under the umbrella (well-chosen).

\paragraph{Lower bounds for sub-Gaussian bandit problems in the case when $\mu^\star$ or the gaps~$\Delta$ are known.}
This framework and the exploitation of this knowledge was first studied by~\citemor{BPR}.
They consider a bandit model $\mathcal{D}$ containing only sub-Gaussian distributions with parameter $\sigma^2 \leq 1$;
that is, distributions $\nu_a$, with expectations $\mu_a \in \R$, such that
\begin{equation}
\label{eq:sub-Gaussian}
\forall \lambda \in \R, \qquad \int_{\R} \exp\bigl( \lambda(y-\mu_a) \bigr) \,
\mathrm{d}\nu_a(y) \leq \exp\biggl(\frac{\lambda^2}{2}\biggr)\,.
\end{equation}
Examples of such distributions include Gaussian distributions with variance smaller than~$1$
and bounded distributions with range smaller than~$2$.

They study how much smaller the regret bounds can get when either the maximal expected payoff
$\mu^\star$ or the gaps $\Delta_a$
are known. For the case when the gaps $\Delta_a$ are known but not $\mu^\star$,
they exhibit a lower bound on the regret matching previously known upper bounds,
thus proving their optimality.
For the case when $\mu^\star$ is known but not the gaps, they
offer an algorithm and its associated regret upper bound,
as well as a framework for deriving a lower
bound; later work (see \citemorpar{err} and \citemorpar{ESAIM})
point out that a bounded regret can be achieved in this case.

We (re-)derive these two lower bounds in a few lines in Section~\ref{sec:rederiv-BPR} of the appendix.
In particular, the well-chosen difficult bandit problems used
are composed of Gaussian distributions $\mathcal{N}(\mu_a,1)$,
with expectations $\mu_a \in \{-\Delta,0,\Delta\}$.
Only statements of the form (well-chosen), not of the form (all), are obtained.
Put differently,
no general distribution-dependent statement like: ``For all bandit problems
in which the gaps $\Delta$ (or the maximal expected payoff $\mu^\star$) are known,
all (reasonable) strategies suffer a regret larger than [...]'' is proposed
by~\citemor{BPR}; only well-chosen, difficult bandit problems are considered.
This is in strong contrast with our general distribution-dependent bounds
for the initial linear regime, provided in Section~\ref{sec:smallT}.

\subsection{Outline of our contributions.}

In Section~\ref{sec:fundamental}, we present Inequality~\eqref{eq:funda}, in our opinion the most efficient and most versatile tool for proving lower bounds in bandit models. We carefully detail its remarkably simple proof, together with an elegant re-derivation of the earlier asymptotic lower bounds by \citemor{LaRo85}, \citemor{BuKa96} and \citemor{CoKa}. Some other earlier bounds are also re-derived in Appendix~\ref{sec:appLB}, namely, the distribution-free lower bound by~\citemor{AuCBFrSc02} as well as the bounded-regret Gaussian lower bounds by~\citemor{BPR} in the case when $\mu^\star$ or the gaps $\Delta$
are known.

The true power of Inequality~\eqref{eq:funda} is illustrated in Section~\ref{sec:smallT}: we study the initial regime when the small number $T$ of draws does not yet permit to unambiguously identify the best arm. We propose three different bounds (each with specific merits). They explain the quasi-linear growth of the regret in this initial phase. We also discuss how the length of the initial phase depends on the number of arms and on the gap between optimal and sub-optimal arms in Kullback-Leibler divergence. These lower bounds are extremely strong as they hold for all possible bandit problems, not
just for some well-chosen ones.

Section~\ref{sec:largeT} contains a general non-asymptotic lower bound for the logarithmic (large $T$) regime. This bound does not only contain the right leading term, but the analysis aims at highlighting what the second-order terms depend on. Results of independent interest on the regularity (upper semi-continuity) of $\Kinf$ are provided
in its Subsection~\ref{sec:wbm}.

\section{The fundamental inequality, and re-derivation of earlier lower bounds.}
\label{sec:fundamental}

We denote by $\kl$ the Kullback-Leibler divergence for
Bernoulli distributions:
\begin{equation}
\label{eq:kl}
\forall p,q \in [0,1]^2, \qquad \kl(p,q) = p \ln \frac{p}{q} + (1-p) \ln \frac{1-p}{1-q}\,.
\end{equation}
We show in this section that for all strategies $\psi$,
for all bandit problems $\unu$ and $\unu'$,
for all $\sigma(I_T)$--measurable random variables $Z$ with values in $[0,1]$,
\begin{equation}
\label{eq:funda}
\sum_{a=1}^K \E_{\unu} \! \bigl[ N_{\psi,a}(T) \bigr] \, \KL(\nu_a,\nu'_a)
\geq \kl\bigl( \E_{\unu}[Z],\,\E_{\unu'}[Z] \bigr)\,.
\end{equation}

Inequality~\eqref{eq:funda} will be referred to as the fundamental inequality
of this article.
We will typically apply it by considering variables of the form
$Z = N_{\psi,k}(T)/T$ for some arm $k$.
That the $\kl$ term in~\eqref{eq:funda} then also contains expected numbers of draws of arms
will be very handy. Unlike all previous proofs of distribution-dependent lower bounds
for bandit problems, we will not have to introduce well-chosen events and control their
probability by some Markov--Chernoff bounding.
Implicit changes of measures will however be performed by
considering bandit problems $\unu$ and $\unu'$ and their
associated probability measures $\mathbb{P}_{\unu}$ and $\mathbb{P}_{\unu'}$.

\paragraph{Underlying probability measures.}
The proof of~\eqref{eq:funda} will be based, among others, on an application
of the chain rule for Kullback-Leibler divergences. For this reason, it is
helpful to construct and define the underlying measures, so that the needed
stochastic transition kernels appear clearly.

By Kolmogorov's extension
theorem, there exists a measurable space $(\Omega,\mathcal{F})$
on which all probability measures $\mathbb{P}_{\unu}$ and $\mathbb{P}_{\unu'}$
considered above can be defined; e.g., $\Omega = [0,1] \times \bigl( \R \times [0,1] \bigr)^\mathbb{N}$.
Given the probabilistic and strategic setting described in Section~\ref{sec:setting},
the probability measure $\P_{\unu}$ over this $(\Omega,\mathcal{F})$ is such that
for all $t \geq 0$,
for all Borel sets $B \subseteq \R$ and $B' \subseteq [0,1]$,
\begin{equation}
\label{eq:repr}
\P_{\unu} \bigl( Y_{t+1} \in B, \, U_{t+1}\in B' \,\, \big| \,\, I_t \bigr)
= \nu_{\psi_t(I_t)}(B) \,\, \lambda(B')\,,
\end{equation}
where $\lambda$ denotes the Lebesgue measure on $[0,1]$.

\begin{rema}
Equation~\eqref{eq:repr} actually reveals that the
distributions $\P_{\unu}$ should be indexed as well
by the considered strategy $\psi$. Because the important element
in the proofs will be the dependency on $\unu$ (we will replace
$\unu$ by alternative bandit problems $\unu'$), we drop the
dependency on $\psi$ in the notation for the underlying probability measures.
This will not come at the cost of clarity as virtually all
events $A_\psi$ and random variables $Z_\psi$ that will be considered will depend on $\psi$:
we will almost always deal with probabilities of the form
$\P_{\unu}(A_\psi)$ or expectations of the form $\E_{\unu}[Z_\psi]$.
\end{rema}

\subsection{Proof of the fundamental inequality~\eqref{eq:funda}.}

We let $\P_{\unu}^{I_T}$ and $\P_{\unu'}^{I_T}$ denote the respective distributions (pushforward measures)
of $I_T$ under $\P_{\unu}$ and $\P_{\unu'}$.
We add an intermediate equation in~\eqref{eq:funda},
\begin{equation}
\label{eq:funda-long}
\sum_{a=1}^K \E_{\unu} \! \bigl[ N_{\psi,a}(T) \bigr] \, \KL(\nu_a,\nu'_a)
= \KL\bigl( \P_{\unu}^{I_T}, \, \P_{\unu'}^{I_T} \bigr)
\geq \kl\bigl( \E_{\unu}[Z],\,\E_{\unu'}[Z] \bigr)\,,
\end{equation}
and are left with proving a standard equality (via the chain rule
for Kullback-Leibler divergences) and a less standard inequality
(following from the data-processing inequality for Kullback-Leibler divergences).

\begin{remark}
Although this possibility is not used in the present article, it is important to note, after~\citemoradd{KaCaGa16}{Lemma 1}, that \eqref{eq:funda-long} actually holds not only for deterministic values of $T$ but also for any stopping time with respect to the filtration generated by $(I_t)_{t \geq 1}$.
\end{remark}

\paragraph{Proof of the equality in~\eqref{eq:funda-long}.}

This equality can be found, e.g., in the proofs of the distribution-free lower bounds on the bandit regret,
in the special case of Bernoulli distributions,
see \citemor{AuCBFrSc02} and \citemor{CBLu06}; see also~\citemor{CoPr14}.
We thus reprove this equality for the sake of completeness only.

We use the symbol $\otimes$ to denote products of measures.
The stochastic transition kernel~\eqref{eq:repr} exactly indicates that
the conditional distribution of $(Y_{t+1},U_{t+1})$ given $I_t$ equals
\[
\P_{\unu}^{(Y_{t+1},U_{t+1})\,|\,I_t} = \nu_{\psi_t(I_t)} \otimes \lambda\,.
\]
Because the conditional distribution at hand takes such a simple form,
the chain rule for Kullback-Leibler divergences applies; it ensures that
for all $t \geq 0$,
\begin{align}
\nonumber
\KL\Bigl( \P_{\unu}^{I_{t+1}}, \, \P_{\unu'}^{I_{t+1}} \Bigr)
& = \KL\Bigl( \P_{\unu}^{(I_t,Y_{t+1},U_{t+1})}, \, \P_{\unu'}^{(I_t,Y_{t+1},U_{t+1})} \Bigr) \\
\label{eq:chainruleiterate}
& = \KL\bigl( \P_{\unu}^{I_t}, \, \P_{\unu'}^{I_t} \bigr)
+ \KL\Bigl( \P_{\unu}^{(Y_{t+1},U_{t+1})\,|\,I_t}, \, \P_{\unu'}^{(Y_{t+1},U_{t+1})\,|\,I_t} \Bigr)\,,
\end{align}
where
\begin{align*}
\KL\Bigl( \P_{\unu}^{(Y_{t+1},U_{t+1})\,|\,I_t}, \, \P_{\unu'}^{(Y_{t+1},U_{t+1})\,|\,I_t} \Bigr)
& = \E_{\unu} \! \biggl[ \E_{\unu} \! \Bigl[ \KL\bigl( \nu_{\psi_t(I_t)} \otimes \lambda, \,
\nu'_{\psi_t(I_t)} \otimes \lambda \bigr) \, \Big| \, I_t \Bigr] \biggr] \\
& = \E_{\unu} \! \biggl[ \E_{\unu} \! \Bigl[ \KL\bigl( \nu_{\psi_t(I_t)}, \,
\nu'_{\psi_t(I_t)} \bigr) \, \Big| \, I_t \Bigr] \biggr] \\
& = \E_{\unu} \!\! \left[ \sum_{a=1}^K \KL(\nu_a,\nu'_a) \, \ind_{\{\psi_t(I_t) = a \}} \right].
\end{align*}
Recalling that $A_{t+1} = \psi_t(I_t)$, we proved so far
\[
\KL\Bigl( \P_{\unu}^{I_{t+1}}, \, \P_{\unu'}^{I_{t+1}} \Bigr) =
\KL\bigl( \P_{\unu}^{I_t}, \, \P_{\unu'}^{I_t} \bigr) +
\E_{\unu} \!\! \left[ \sum_{a=1}^K \KL(\nu_a,\nu'_a) \, \ind_{\{A_{t+1} = a \}} \right].
\]
Iterating the argument
and using that
$\KL\bigl( \P_{\unu}^{I_0}, \, \P_{\unu'}^{I_0} \bigr) =
\KL\bigl( \P_{\unu}^{U_0}, \, \P_{\unu'}^{U_0} \bigr) =
\KL(\lambda,\lambda)=0$
leads to the equality stated in~\eqref{eq:funda-long}.

\paragraph{Proof of the inequality in~\eqref{eq:funda-long}.}

\emph{This is our key contribution to a simplified proof of the lower bound}~\eqref{eq:obj}.
It is a consequence of the data-processing inequality (also known as contraction of entropy),
i.e., the fact that Kullback-Leibler
divergences between pushforward measures are smaller than the Kullback-Leibler
divergences between the original probability measures; see Lemma~\ref{lm:DPineq} in
Appendix~\ref{sec:reminderinfotheory} for a statement and elements of proof.

We actually state our inequality in a slightly more general way, as it is of independent interest.

\begin{lemma}
\label{lm:DPE}
Consider a measurable space $(\Gamma,\mathcal{G})$ equipped with two distributions
$\P_1$ and $\P_2$, and any $\mathcal{G}$--measurable random variable $Z : \Omega \to [0,1]$.
We denote respectively by $\E_1$ and $\E_2$ the expectations under $\P_1$ and $\P_2$.
Then,
\[
\KL(\P_1,\P_2) \geq \kl\bigl(\E_1[Z],\E_2[Z]\bigr)\,.
\]
\end{lemma}

\begin{MORproof}
We augment the underlying measurable space into $\Gamma \times [0,1]$, where
$[0,1]$ is equipped with the Borel $\sigma$--algebra $\cB\bigl([0,1]\bigr)$ and the Lebesgue measure $\lambda$.
We denote by $\cG \otimes \cB\bigl([0,1]\bigr)$ the $\sigma$--algebra generated by product sets in $\cG \times \cB\bigl([0,1]\bigr)$.
Now, for any event $E \in \cG \otimes \cB\bigl([0,1]\bigr)$,
by the consideration of product distributions for the equality
and by the data-processing inequality (Lemma~\ref{lm:DPineq}) applied to $X = \ind_E$ for the inequality,
we have
\[
\KL(\P_1,\P_2) = \KL\bigl(\P_1 \otimes \lambda,\,\P_2 \otimes \lambda \bigr)
\geq \KL\Bigl( (\P_1 \otimes \lambda)^{\ind_E},\,\,(\P_2 \otimes \lambda)^{\ind_E} \Bigr)\,.
\]
The distribution $(\P_j \otimes \lambda)^{\ind_E}$ of $\ind_E$ under $\P_j \otimes \lambda$ is a Bernoulli
distribution, with parameter the probability of $E$ under $\P_j \otimes \lambda$; therefore,
using the notation $\kl$, we have got so far
\[
\KL(\P_1,\P_2) \geq \KL\Bigl( (\P_1 \otimes \lambda)^{\ind_E},\,\,(\P_2 \otimes \lambda)^{\ind_E} \Bigr)
= \kl \bigl( (\P_1 \otimes \lambda)(E), \, (\P_2 \otimes \lambda)(E) \bigr)\,.
\]
We consider $E = \bigl\{ (\gamma,x) \in \Gamma \times [0,1] : x \leq Z(\gamma) \bigr\}$ and note
noting that for all $j$, by the Fubini-Tonelli theorem,
\[
(\P_j \otimes \lambda)(E) =
\bigintsss_\Omega \left( \int_{[0,1]} \ind_{\{ x \leq Z(\gamma) \}} \,\d\lambda(x) \right) \d\P_j(\gamma)
= \int_\Omega Z(\gamma) \,\d\P_j(\gamma)
= \E_j[Z]\,.
\]
This concludes the proof of this lemma.
\end{MORproof}

\subsection{Application: re-derivation of the general asymptotic distribution-dependent bound.}
\label{sec:rederiv-asympt}

As a warm-up, we show how the asymptotic distribution-dependent lower bound~\eqref{eq:obj}
of \citemor{BuKa96} can be reobtained, for so-called uniformly fast convergent strategies.

\begin{defi}
\label{eq:consistency}
A strategy $\psi$ is uniformly fast convergent on a model $\cD$ if for all bandit problems $\unu$ in $\cD$, for all
suboptimal arms $a$, i.e., for all arms $a$ such that $\Delta_a > 0$, for all $0 < \alpha \leq 1$,
it satisfies $\E_{\unu}\bigl[N_{\psi,a}(T)\bigr] = o(T^\alpha)$.
\end{defi}

\begin{theorem}
\label{th:asympt}
For all models $\mathcal{D}$, for all uniformly fast convergent strategies $\psi$ on $\cD$,
for all bandit problems $\unu$, for all suboptimal arms $a$,
\[
\liminf_{T \to \infty} \,\, \frac{\E_{\unu} \! \bigl[ N_{\psi,a}(T) \bigr]}{\ln T}
\geq \frac{1}{\Kinf(\nu_a,\mu^\star,\cD)}\,.
\]
\end{theorem}

\begin{MORproof}
Given any bandit problem $\unu$ and any suboptimal arm $a$, we consider a modified problem $\unu'$
where $a$ is the (unique) optimal arm: $\nu'_k = \nu_k$ for all $k \ne a$ and
$\nu'_a$ is any distribution in $\mathcal{D}$ such that its expectation $\mu'_a$ satisfies
$\mu'_a > \mu^\star$ (if such a distribution exists; see the end of the proof otherwise).
We apply the fundamental inequality~\eqref{eq:funda} with $Z = N_{\psi,a}(T)/T$.
All Kullback-Leibler divergences in its
left-hand side are null except the one for arm $a$, so that we get the lower bound
\begin{align}
\nonumber
\E_{\unu} \! \bigl[ N_{\psi,a}(T) \bigr] \, \KL(\nu_a,\nu'_a)
& \geq \kl\Biggl( \frac{\E_{\unu} \! \bigl[ N_{\psi,a}(T) \bigr]}{T}, \,\, \frac{\E_{\unu'} \! \bigl[ N_{\psi,a}(T) \bigr]}{T} \Biggr) \\
\label{eq:BuKa-almost2}
& \geq \left( 1 - \frac{\E_{\unu} \! \bigl[ N_{\psi,a}(T) \bigr]}{T} \right)
\ln \frac{T}{T - \E_{\unu'} \! \bigl[ N_{\psi,a}(T) \bigr]} - \ln 2\,,
\end{align}
where we used for the second inequality that for all $(p,q) \in [0,1]^2$,
\begin{equation}
\kl(p,q) = \underbrace{p \ln \frac{1}{q}}_{\geq 0} + (1-p) \ln \frac{1}{1-q} +
\bigl( \underbrace{p \ln p + (1-p) \ln(1-p)}_{\geq - \ln 2} \bigr)\,.
\label{ineg:kl}
\end{equation}
The uniform fast convergence of $\psi$ together with
the fact that all arms $k \ne a$ are suboptimal for $\unu'$
entails that
\[
\forall \, 0 < \alpha \leq 1, \qquad
0 \leq T - \E_{\unu'} \! \bigl[ N_{\psi,a}(T) \bigr] = \sum_{k \ne a}
\E_{\unu'} \! \bigl[ N_{\psi,k}(T) \bigr] = o(T^\alpha)\,;
\]
in particular,
$T - \E_{\unu'} \! \bigl[ N_{\psi,a}(T) \bigr] \leq T^\alpha$ for $T$ sufficiently large.
Therefore, for all $0 < \alpha \leq 1$,
\[
\liminf_{T \to \infty} \,\, \frac{1}{\ln T} \ln \frac{T}{T - \E_{\unu'} \! \bigl[ N_{\psi,a}(T) \bigr]}
\geq \liminf_{T \to \infty} \,\, \frac{1}{\ln T} \ln \frac{T}{T^\alpha} = (1-\alpha)\,.
\]
In addition, the uniform fast convergence of $\psi$
and the suboptimality of $a$ for the bandit problem $\unu$ ensure that
$\E_{\unu} \! \bigl[ N_{\psi,a}(T) \bigr] / T \to 0$. Substituting these two facts
in~\eqref{eq:BuKa-almost2} we proved
\[
\liminf_{T \to \infty} \,\, \frac{\E_{\unu} \! \bigl[ N_{\psi,a}(T) \bigr]}{\ln T}
\geq \frac{1}{\KL(\nu_a,\nu'_a)}\,.
\]
By taking the supremum in the right-hand side
over all distributions $\nu'_a \in \mathcal{D}$ with $\mu'_a > \mu^\star$, if at least one
such distribution exists, we get the bound of the theorem.
Otherwise, $\Kinf(\nu_a,\mu^\star,\cD) = +\infty$ by a standard
convention on the infimum of an empty set and the bound holds as well.
\end{MORproof}

\section{Non-asymptotic bounds for small values of $T$.}
\label{sec:smallT}

We prove three such bounds with different merits and drawbacks.
Basically, we expect suboptimal arms to be pulled each about $T/K$
of the time when $T$ is small; when $T$ becomes larger, sufficient information was gained for identifying the best arm, and the logarithmic regime can take place.

The first bound shows that $\E_{\unu} \! \bigl[ N_{\psi,a}(T) \bigr]$ is of order $T/K$  as long as $T$ is at most of order $1/\Kinf(\nu_a,\mu^\star,\cD)$;
we call it an absolute lower bound for a suboptimal arm $a$.
Its drawback is that the times $T$ for which it is valid are independent of the number of arms $K$, while (at least in some cases) one may expect the initial phase to last until $T \approx K/\Kinf(\nu_a,\mu^\star,\cD)$.

The second lower bound thus addresses the dependency of the initial phase in $K$ by considering a relative lower bound between a suboptimal arm $a$ and an optimal arm $a^\star$. We prove that $\E_{\unu} \! \bigl[ N_{\psi,a}(T)/N_{\psi,a^\star}(T) \bigr]$ is not much smaller than $1$ whenever $T$ is at most of order $K/\KL(\nu_a,\nu_{a^\star})$. Here, the number of arms $K$ plays the expected effect on the length of the initial exploration phase, which
should be proportional to $K$.

The third lower bound is a collective lower bound on all suboptimal arms,
i.e., a lower bound on $\sum_{a \not\in \cA^\star(\unu)} \E_{\unu} \! \bigl[N_{\psi,a}(T)\bigr]$
where $\cA^\star(\unu)$ denotes the set of the $A^\star_{\unu}$ optimal arms of $\unu$.
It is of the desired order $T(1-A^\star_{\unu}/K)$ for times $T$
of the desired order $K/\cK^{\max}_{\unu}$, where $\cK^{\max}_{\unu}$ is some Kullback-Leibler divergence.

\paragraph{Minimal restrictions on the considered strategies.}
We prove these lower bounds under minimal assumptions on the considered strategies:
either some mild symmetry (much milder than asking for symmetry under permutation
of the arms, see Definition~\ref{def:ps});
or the fact that for suboptimal arms $a$, the number of pulls
$\E_{\unu} \bigl[ N_{\psi,a}(T) \bigr]$ should decrease as $\mu_a$ decreases,
all other distributions of arms being fixed (see Definitions~\ref{def:smarter}
and~\ref{def:monotonic}).
These assumptions are satisfied by all well-performing strategies
we could think of: the UCB strategy of \citemor{AuCBFi02},
the KL-UCB strategy of \citemor{klucb},
\citemor{Th33} Sampling, EXP3 of \citemor{AuCBFrSc02}, etc.

These mild restrictions on the considered strategies are necessary to rule out the irrelevant
strategies (e.g., always pull arm~$1$) that would perform extremely well for some
particular bandit problems $\unu$. This is because we aim at proving distribution-dependent lower bounds
that are valid for all bandit problems $\unu$: we prefer to impose the (mild) constraints on the strategies.

Note that the assumption of uniform fast convergence (Definition~\ref{eq:consistency}), though classical and
well accepted, is quite strong. Note that it is necessary for a strategy to satisfy some symmetry
and to be smarter than the uniform strategy in the limit (not for all $T$, see
Definition~\ref{def:smarter}) to be uniformly fast convergent.
Hence, the class of strategies we consider is essentially much larger than the subset of
uniformly fast convergent strategies.

\subsection{Absolute lower bound for a suboptimal arm.}

The uniform strategy is the one that pulls an arm uniformly
at random at each round.

\begin{defi}
\label{def:smarter}
A strategy $\psi$ is smarter than the uniform strategy on a model $\cD$ if
for all bandit problems $\unu$ in $\cD$, for all optimal arms $a^\star$,
for all $T \geq 1$,
\[
\E_{\unu} \bigl[ N_{\psi,a^\star}(T) \bigr] \geq \frac{T}{K}\,.
\]
\end{defi}

\begin{theorem}
\label{th:2}
For all models $\cD$,
for all strategies $\psi$ that are smarter than the uniform strategy on $\cD$,
for all bandit problems~$\unu$, for all arms $a$, for all $T \geq 1$,
\[
\E_{\unu} \bigl[ N_{\psi,a}(T) \bigr] \geq \frac{T}{K} \Bigl( 1 - \sqrt{2 T \Kinf(\nu_a,\mu^\star,\cD)} \Bigr)\,.
\]
In particular,
\[
\forall \, T \leq \frac{1}{8 \Kinf(\nu_a,\mu^\star,\cD)}, \quad \qquad
\E_{\unu} \bigl[ N_{\psi,a}(T) \bigr] \geq \frac{T}{2K}\,.
\]
\end{theorem}

\begin{MORproof}
The definition of being smarter than the uniform strategy
takes care of the lower bound for optimal arms $a$:
it thus suffices to consider suboptimal arms $a$.
As in the proof of Theorem~\ref{th:asympt},
we consider a modified bandit problem $\unu'$
with $\nu'_k = \nu_k$ for all $k \ne a$ and
$\nu'_a \in \mathcal{D}$ such that $\mu'_a > \mu^\star$,
if such a distribution $\nu'_a$ exists (otherwise, the first claimed lower
bounds equals $-\infty$).
From~\eqref{eq:funda}, we get
\[
\E_{\unu} \! \bigl[ N_{\psi,a}(T) \bigr] \, \KL(\nu_a,\nu'_a)
\geq \kl\Biggl( \frac{\E_{\unu} \! \bigl[ N_{\psi,a}(T) \bigr]}{T}, \,\, \frac{\E_{\unu'} \! \bigl[ N_{\psi,a}(T) \bigr]}{T} \Biggr)\,.
\]
We may assume that $\E_{\unu} \! \bigl[ N_{\psi,a}(T) \bigr]/T \leq 1/K$; otherwise, the first claimed bound holds.
Since $a$ is the optimal arm under $\unu'$ and since the considered strategy is smarter than the uniform
strategy, $\E_{\unu'} \! \bigl[ N_{\psi,a}(T) \bigr]/T \geq 1/K $.
Using that $q \mapsto \kl(p,q)$ is increasing on $[p,1]$, we thus get
\[
\kl\Biggl( \frac{\E_{\unu} \! \bigl[ N_{\psi,a}(T) \bigr]}{T}, \,\, \frac{\E_{\unu'} \! \bigl[ N_{\psi,a}(T) \bigr]}{T} \Biggr)
\geq \kl\Biggl( \frac{\E_{\unu} \! \bigl[ N_{\psi,a}(T) \bigr]}{T}, \,\, \frac{1}{K} \Biggr)\,.
\]
Lemma~\ref{lm:Pinskerlocal} of Appendix~\ref{sec:reminderinfotheory} yields
\[
\E_{\unu} \! \bigl[ N_{\psi,a}(T) \bigr] \, \KL(\nu_a,\nu'_a) \geq
\kl\Biggl( \frac{\E_{\unu} \! \bigl[ N_{\psi,a}(T) \bigr]}{T}, \,\, \frac{1}{K} \Biggr) \geq
\frac{K}{2} \Biggl( \frac{\E_{\unu} \! \bigl[ N_{\psi,a}(T) \bigr]}{T} - \frac{1}{K} \Biggr)^{\!\! 2},
\]
from which follows, after substitution of the above assumption $\E_{\unu} \! \bigl[ N_{\psi,a}(T) \bigr]/T \leq 1/K$
in the left-hand side,
\[
\frac{\E_{\unu} \! \bigl[ N_{\psi,a}(T) \bigr]}{T} \geq
\frac{1}{K} - \sqrt{\frac{2T}{K^2} \,\KL(\nu_a,\nu'_a)}\,.
\]
Taking the supremum of the right-hand side over all $\nu'_a \in \cD$ such that $E(\nu'_a) > \mu^\star$ and rearranging
concludes the proof.
\end{MORproof}

\subsection{Relative lower bound.}

Our proof will be based on an assumption of symmetry (milder than
requiring that if the arms are permuted in a bandit problem,
the algorithm behaves the same way, as in Definition~\ref{def:invtrans:sym}).

\begin{defi}
\label{def:ps}
A strategy $\psi$ is pairwise symmetric for optimal arms on $\cD$
if for all bandit problems $\unu$ in $\cD$, for each pair of
optimal arms $a^\star$ and $a_\star$, the equality
$\nu_{a^\star} = \nu_{a_\star}$ entails that, for all $T \geq 1$,
\[
\bigl( N_{\psi,a^\star}(T), \, N_{\psi,a_\star}(T) \bigr)
\qquad \mbox{and} \qquad
\bigl( N_{\psi,a_\star}(T), \, N_{\psi,a^\star}(T) \bigr)
\]
have the same distribution.
\end{defi}

Note that the required symmetry is extremely mild as only pairs of \emph{optimal} arms
with the \emph{same} distribution are to be considered. What the equality of distributions
means is that the strategy should be based only on payoffs and not on the values of the
indexes of the arms.

\begin{theorem}
\label{th:3}
For all models $\cD$,
for all strategies $\psi$ that are pairwise symmetric for optimal arms on $\cD$,
for all bandit problems $\unu$ in $\cD$, for all suboptimal arms $a$ and all optimal arms $a^\star$,
for all $T \geq 1$,
\[
\mbox{either} \quad \E_{\unu} \bigl[ N_{\psi,a}(T) \bigr] \geq \frac{T}{K}
\qquad
\mbox{or} \quad \E_{\unu} \!\!\left[
\frac{\max\bigl\{N_{\psi,a}(T),\,1\bigr\}}{\max\bigl\{N_{\psi,a^\star}(T),\,1\bigr\}} \right]
\geq 1 - 2 \sqrt{\frac{2T \, \KL(\nu_a,\nu_{a^\star})}{K}} \,.
\]
In particular,
\[
\forall \, T \leq \frac{K}{32 \, \KL(\nu_a,\nu_{a^\star})}, \quad \qquad
\mbox{either} \quad \E_{\unu} \bigl[ N_{\psi,a}(T) \bigr] \geq \frac{T}{K}
\qquad
\mbox{or} \quad \E_{\unu} \!\!\left[
\frac{\max\bigl\{N_{\psi,a}(T),\,1\bigr\}}{\max\bigl\{N_{\psi,a^\star}(T),\,1\bigr\}} \right]
\geq \frac{1}{2}\,.
\]
\end{theorem}

That is, on average, in the small $T$ regime, each suboptimal arm is played at least half the number of
times when an optimal arm was played.

\begin{MORproof}
For all arms $k$, we denote by $N^+_{\psi,k}(T) = \max\bigl\{N_{\psi,k}(T),\,1\bigr\}$.
Given a bandit problem $\unu$ and a suboptimal arm $a$, we form an alternative bandit problem
$\unu'$ given by $\nu'_k = \nu_k$ for all $k \ne a$ and $\nu'_a = \nu_{a^\star}$,
where $a^\star$ is an optimal arm of $\unu$. In particular, arms $a$ and $a^\star$ are both optimal arms
under $\unu'$. By the assumption of pairwise symmetry for optimal arms, we have in particular that
\[
\E_{\unu'} \!\!\left[ \frac{N^+_{\psi,a}(T)}{N^+_{\psi,a}(T)+N^+_{\psi,a^\star}(T)} \right]
= \E_{\unu'} \!\!\left[ \frac{N^+_{\psi,a^\star}(T)}{N^+_{\psi,a^\star}(T)+N^+_{\psi,a}(T)} \right]
= \frac{1}{2}\,.
\]
The latter equality and
the fundamental inequality~\eqref{eq:funda} yield in the present case,
through the choice of $Z = N^+_{\psi,a}(T)\big/\bigl(N^+_{\psi,a}(T)+N^+_{\psi,a^\star}(T)\bigr)$,
\begin{equation}
\label{eq:funda-Tsmall2}
\E_{\unu} \! \bigl[ N_{\psi,a}(T) \bigr] \,\, \KL(\nu_a,\nu'_a)
\geq \kl\!\left( \E_{\unu} \!\!\left[ \frac{N^+_{\psi,a}(T)}{N^+_{\psi,a}(T)+N^+_{\psi,a^\star}(T)} \right], \,\,
\frac{1}{2} \right).
\end{equation}
The concavity of the function $x \mapsto x/(1+x)$ and Jensen's inequality show that
\[
\E_{\unu} \!\!\left[ \frac{N^+_{\psi,a}(T)}{N^+_{\psi,a}(T)+N^+_{\psi,a^\star}(T)} \right]
= \E_{\unu} \!\!\left[ \frac{N^+_{\psi,a}(T)\big/N^+_{\psi,a^\star}(T)}{1+N^+_{\psi,a}(T)\big/N^+_{\psi,a^\star}(T)} \right]
\leq \frac{ \E_{\unu}\bigl[N^+_{\psi,a}(T)\big/N^+_{\psi,a^\star}(T)\bigr] }{1 + \E_{\unu}\bigl[N^+_{\psi,a}(T)\big/N^+_{\psi,a^\star}(T)\bigr]}\,.
\]
We can assume that $\E_{\unu} \! \bigl[N^+_{\psi,a}(T)\big/N^+_{\psi,a^\star}(T)\bigr] \leq 1$, otherwise, the result
of the theorem is obtained. In this case, the latter upper bound is smaller than $1/2$.
Using in addition that $p \mapsto \kl(p,1/2)$ is decreasing on $[0,1/2]$,
and assuming that $\E_{\unu} \! \bigl[ N_{\psi,a}(T) \bigr] \leq T/K$ (otherwise, the result of the theorem
is obtained as well), we get from~\eqref{eq:funda-Tsmall2}
\[
\frac{T}{K} \,\, \KL(\nu_a,\nu'_a)
\geq \kl\!\left( \frac{ \E_{\unu} \! \bigl[N^+_{\psi,a}(T)\big/N^+_{\psi,a^\star}(T)\bigr] }{1 +
\E_{\unu} \! \bigl[N^+_{\psi,a}(T)\big/N^+_{\psi,a^\star}(T)\bigr]}, \,\,
\frac{1}{2} \right).
\]
Pinsker's inequality (in its classical form, see Appendix~\ref{sec:reminderinfotheory} for a statement) entails the inequality
\[
\frac{T}{K} \,\, \KL(\nu_a,\nu'_a)
\geq 2 \left( \frac{1}{2} - \frac{r}{1+r} \right)^{\!\! 2}, \qquad
\mbox{where} \quad r = \E_{\unu} \!\!\left[ \frac{N^+_{\psi,a}(T)}{N^+_{\psi,a^\star}(T)} \right]\,.
\]
In particular,
\[
\frac{r}{1+r} \geq \frac{1}{2} - \sqrt{\frac{T \, \KL(\nu_a,\nu'_a)}{2 K}}\,.
\]
Applying the increasing function $x \mapsto x/(1-x)$ to both sides, we get
\[
r \geq \frac{1 - \sqrt{2T \, \KL(\nu_a,\nu'_a)/K}}{1 + \sqrt{2T \, \KL(\nu_a,\nu'_a)/K}}
\geq \left( 1 - \sqrt{\frac{2T \, \KL(\nu_a,\nu'_a)}{K}} \right)^{\!\! 2},
\]
where we used $1/(1+x) \geq 1-x$ for the last inequality and where we assumed
that $T$ is small enough to ensure $1 - \sqrt{2T \, \KL(\nu_a,\nu'_a)/K} \geq 0$.
Whether this condition is satisfied or not, we have the (possibly void) lower bound
\[
r \geq 1 - 2 \sqrt{\frac{2T \, \KL(\nu_a,\nu'_a)}{K}}\,.
\]
The proof is concluded by noting that by definition $\nu'_a = \nu_{a^\star}$.
\end{MORproof}

\subsection{Collective lower bound.}

In this section, for any given bandit problem $\unu$, we denote by
$\cA^\star(\unu)$ the set of its optimal arms and by
$\cW(\unu)$ the set of its worst arms, i.e., the ones associated with the distributions with the
smallest expectation among all distributions for the arms.
We also let $A^\star_{\unu}$ be the cardinality of $\cA^\star(\unu)$.

We define the following partial order $\preccurlyeq$ on bandit problems:
$\unu' \preccurlyeq \unu$ if
\[
\forall a \in \cA^\star(\unu), \quad \nu_a = \nu'_a
\qquad \mbox{and} \qquad
\forall a \not\in \cA^\star(\unu), \quad E(\nu'_a) \leq E(\nu_a)\,.
\]
In particular, $\cA^\star(\unu) = \cA^\star(\unu')$ in this case.
The definition models the fact that the bandit problem $\unu'$ should be easier
than $\unu$, as non-optimal arms in $\unu'$ are farther away from the optimal arms
(in expectation) that in $\unu$. Any reasonable strategy should perform
better on $\unu'$ than on $\unu$, which leads to the following definition,
where we measure performance in the expected number of times optimal arms are pulled.
(Recall that the sets of optimal arms are identical for $\unu$ and $\unu'$.)

\begin{defi}
\label{def:monotonic}
A strategy $\psi$ is monotonic on a model $\cD$ if for all bandit problems $\unu' \preccurlyeq \unu$
in $\cD$,
\[
\sum_{a^\star \in \cA^\star(\unu')} \E_{\unu'}\bigl[N_{\psi,a^\star}(T)\bigr] \geq
\sum_{a^\star \in \cA^\star(\unu)} \E_{\unu}\bigl[N_{\psi,a^\star}(T)\bigr]\,.
\]
\end{defi}

\begin{theorem}
\label{th:4}
For all models $\cD$,
for all strategies $\psi$ that are pairwise symmetric for optimal arms and monotonic on $\cD$, for all bandit problems $\unu$ in $\cD$,
suboptimal arms are collectively sampled at least
\begin{multline*}
\sum_{a \not\in \cA^\star(\unu)} \E_{\unu}\bigl[N_{\psi,a}(T)\bigr]\geq T \Bigg(1-\frac{A^\star_{\unu}}{K} -\frac{A^\star_{\unu} \sqrt{ 2T \, \cK^{\max}_{\unu}}}{K}-\frac{2 A^\star_{\unu} T \cK^{\max}_{\unu}}{K}\Bigg)\,, \\
\mbox{where}
\qquad
\cK^{\max}_{\unu} = \min_{w \in \cW(\unu)} \max_{a^\star \in \cA^\star(\unu)} \KL(\nu_{w},\nu_{a^\star})\,.
\end{multline*}
In particular,
\[
\forall \, T \leq \frac{K}{8 \, A^\star_{\unu} \, \cK^{\max}_{\unu}}, \quad \qquad
\sum_{a \not\in \cA^\star(\unu)} \E_{\unu}\bigl[N_{\psi,a}(T)\bigr]\geq \frac{T}{2} \Bigg(1-\frac{A^\star_{\unu}}{K} \Biggr)\,.
\]
\end{theorem}

To get a lower bound on the regret from this theorem,
we use
\begin{equation}
\label{eq:LBth4}
R_{\psi,\unu,T} \geq \left( \min_{a \not\in \cA^\star(\unu)} \Delta_a \right)
\sum_{a \not\in \cA^\star(\unu)} \E_{\unu}\bigl[N_{\psi,a}(T)\bigr]
\,.
\end{equation}

\begin{MORproof}
We denote by $\tw$ some $w \in \cW(\unu)$ achieving the minimum in the defining equation of $\cK^{\max}_{\unu}$.
We construct two bandit models from $\unu$.
First, the model $\nuun$ differs from $\unu$ only at suboptimal arms $a \not\in \cA^\star(\unu)$, which we associate with
$\nuun_a = \nu_{\tw}$. By construction, $\nuun \preccurlyeq \unu$.

In the second model $\nuzero$, each arm is associated with $\nu_{\tw}$, i.e., $\nuzero_a=\nu_{\tw}$ for all
$a \in \{1,\ldots,K\}$.

By monotonicity of $\psi$,
\[
\sum_{a \not\in \cA^\star(\unu)} \E_{\unu}\bigl[N_{\psi,a}(T)\bigr] \geq
\sum_{a \not\in \cA^\star(\nuun)} \E_{\nuun}\bigl[N_{\psi,a}(T)\bigr]\,.
\]
We can therefore focus our attention, for the rest of the proof, on the $\E_{\nuun}\bigl[N_{\psi,a}(T)\bigr]$.
The strategy is also pairwise symmetric for optimal arms and all arms of $\nuzero$ are optimal.
This implies in particular that $\E_{\nuzero}\bigl[N_{\psi,1}(T)\bigr]= \E_{\nuzero}\bigl[N_{\psi,a}(T)\bigr]$ for all arms $a$,
thus $\E_{\nuzero}\bigl[N_{\psi,a}(T)\bigr]=T/K$ for all arms $a$.

Now, the bound~\eqref{eq:funda} with $Z = \displaystyle{\sum_{a^\star \in \cA^\star(\unu)} \frac{N_{\psi,a^\star}(T)}{T}}$
and the bandit models $\nuzero$ and $\nuun$ gives
\begin{align*}
\sum_{a^\star \in \cA^\star(\unu)} \E_{\nuzero}\bigl[N_{\psi,a^\star}(T)\bigr] \KL(\nu_{\tw},\nu_{a^\star}) & \geq \kl\Bigg(\sum_{a^\star \in \cA^\star(\unu)} \E_{\nuzero}\bigl[N_{\psi,a^\star}(T)\bigr]/T,\sum_{a^\star \in \cA^\star(\unu)} \E_{\nuun}\bigl[N_{\psi,a^\star}(T)\bigr]/T \Bigg)\\
& = \kl \Bigg(\frac{A^\star_{\unu}}{K}, \sum_{a^\star\in \cA^\star(\unu)} \E_{\nuun}\bigl[N_{\psi,a^\star}(T)\bigr]/T\Bigg).
\end{align*}
By definition of $\cK^{\max}_{\unu}$ and $\tw$, and because $\E_{\nuzero}\bigl[N_{\psi,a}(T)\bigr]=T/K$, we have
\[
\sum_{a^\star \in \cA^\star(\unu)} \E_{\nuzero}\bigl[N_{\psi,a^\star}(T)\bigr] \KL(\nu_{\tw},\nu_{a^\star})
\leq \frac{T A^\star_{\unu} \cK^{\max}_{\unu}}{K}\,,
\]
which yields the inequality
\[
\frac{T A^\star_{\unu} \cK^{\max}_{\unu}}{K} \geq \kl\!\left(\frac{A^\star_{\unu}}{K}, \, x\right)
\qquad \mbox{where} \qquad x = \frac{1}{T} \sum_{a^\star\in \cA^\star(\unu)} \E_{\nuun}\bigl[N_{\psi,a^\star}(T)\bigr]\,.
\]
We want to upper bound $x$, in order to get a lower bound on $1-x$.
We assume that $x \geq A^\star_{\unu}/K$, otherwise, the bound~\eqref{eq:bornex} stated below is also satisfied.
Pinsker's inequality (actually, its local refinement stated as Lemma~\ref{lm:Pinskerlocal} in Appendix~\ref{sec:reminderinfotheory}) then ensures that
\[
\frac{T A^\star_{\unu} \cK^{\max}_{\unu}}{K} \geq \frac{1}{2 x}\left(\frac{A^\star_{\unu}}{K}-x\right)^{\!\! 2},
\]
Lemma~\ref{lm:inequation2} below finally entails that
\begin{equation}
\label{eq:bornex}
x \leq \frac{A^\star_{\unu}}{K} \Bigl(1 + 2T \cK^{\max}_{\unu} +\sqrt{2 T \cK^{\max}_{\unu}} \Bigr)\,.
\end{equation}
The proof is concluded by putting all elements together thanks to
the monotonicity of $\psi$ and the definition of $x$:
\[
\sum_{a \not\in \cA^\star(\unu)} \E_{\unu} \! \bigl[N_{\psi,a}(T)\bigr] \geq
\sum_{a \not\in \cA^\star(\nuun)} \E_{\nuun} \! \bigl[N_{\psi,a}(T)\bigr]
= T(1-x)\,.
\]
\vspace{-1cm}

\end{MORproof}

\begin{lemma}
\label{lm:inequation2}
If $x \in \R$ satisfies $(x-\alpha)^2 \leq \beta x$ for some $\alpha \geq 0$
and $\beta \geq 0$, then $x \leq \alpha + \beta + \sqrt{\alpha \beta}$.
\end{lemma}

\begin{MORproof}
By assumption, $x^2 - (2\alpha+\beta)x + \alpha^2 \leq 0$.
We have that $x$ is smaller than the larger root of the associated polynomial, that is,
\[
x \leq \frac{2\alpha + \beta + \sqrt{(2\alpha+\beta)^2 - 4 \alpha^2}}{2}
= \frac{2\alpha + \beta + \sqrt{4\alpha \beta + \beta^2}}{2}\,.
\]
We conclude with
$\sqrt{4\alpha \beta + \beta^2} \leq \sqrt{4\alpha \beta}
+ \sqrt{\beta^2}$.
\end{MORproof}

\subsection{Numerical illustrations.}
In this section we illustrate some of the bounds stated above for the initial linear regime,
namely, the bounds of Theorems~\ref{th:2} and~\ref{th:4}. It turned out that because of the
``or'' statement in Theorem~\ref{th:3}, its bound was less easy to illustrate.
We need much more difficult bandit problems than the one of Figure~\ref{fig:draws} in order to clearly observe the initial linear phase.

\begin{figure}[h!]
\center
\begin{tabular}{lr}
\includegraphics[width=0.47\textwidth]{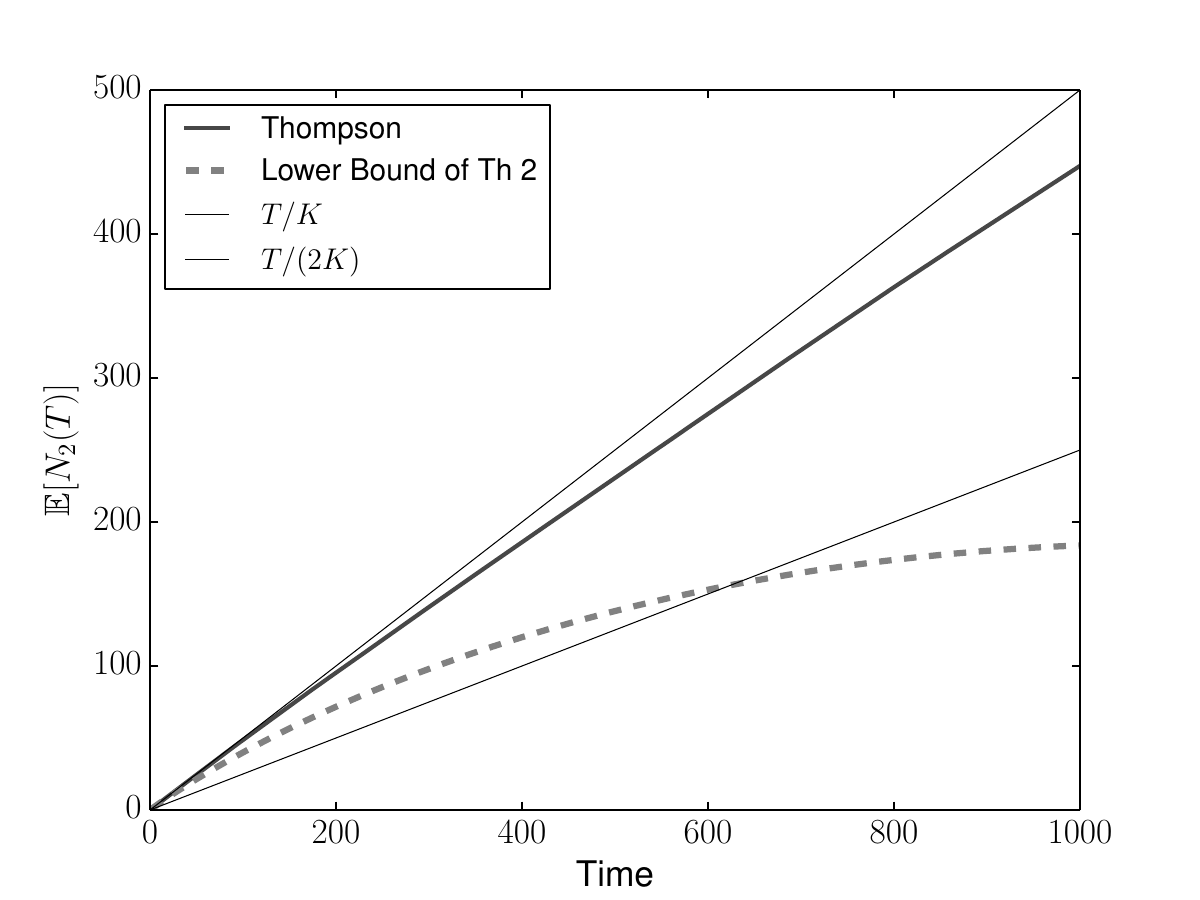} &
\includegraphics[width=0.47\textwidth]{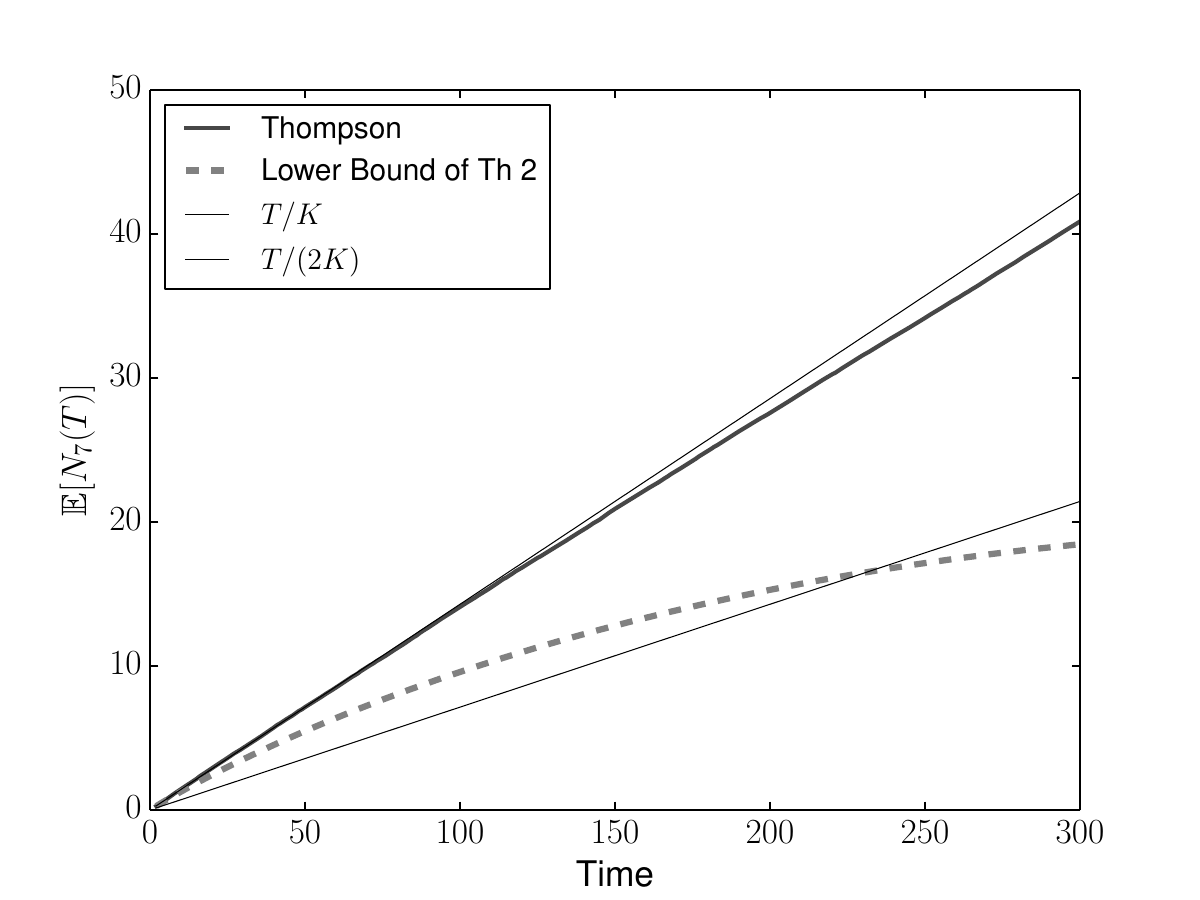}
\end{tabular}
\caption{\label{fig:th2}
Expected number of pulls of the most suboptimal arm for \citemor{Th33} Sampling (\emph{blue, solid} line) on Bernoulli bandit problems, versus the lower bound (\emph{red, dashed} line) of Theorem~\ref{th:2}
for the model $\cD$ of all Bernoulli distributions;
expectations are approximated over $1,000$ runs. \smallskip \newline
\emph{Left}: parameters $(\mu_a)_{1 \leq a \leq 2} =(0.5, \, 0.49)$, with
characteristic time $1/ \big(8 \, \Kinf(\nu_2,\mu^\star,\cD)\big)\approx 625$. \newline
\emph{Right}: parameters $(\mu_a)_{1 \leq a \leq 7} =
(0.05, \, 0.048,\, 0.047, \, 0.046, \, 0.045, \, 0.044, \, 0.043)$, with
$1/ \big(8 \, \Kinf(\nu_7,\mu^\star,\cD) \big)\approx 231$.}
\end{figure}

Theorem~\ref{th:2} is illustrated in Figure~\ref{fig:th2}. We observe that
in the bandit problems contemplated therein, the expected numbers of pulls of the suboptimal
arms considered indeed lie between $T/(2K)$ and $T/K$ in the initial phase, as prescribed
by the theorem. We see, however, that this initial phase is probably longer than what was quantified.

\begin{figure}[!h]
\center
\begin{tabular}{lr}
\includegraphics[width=0.47\textwidth]{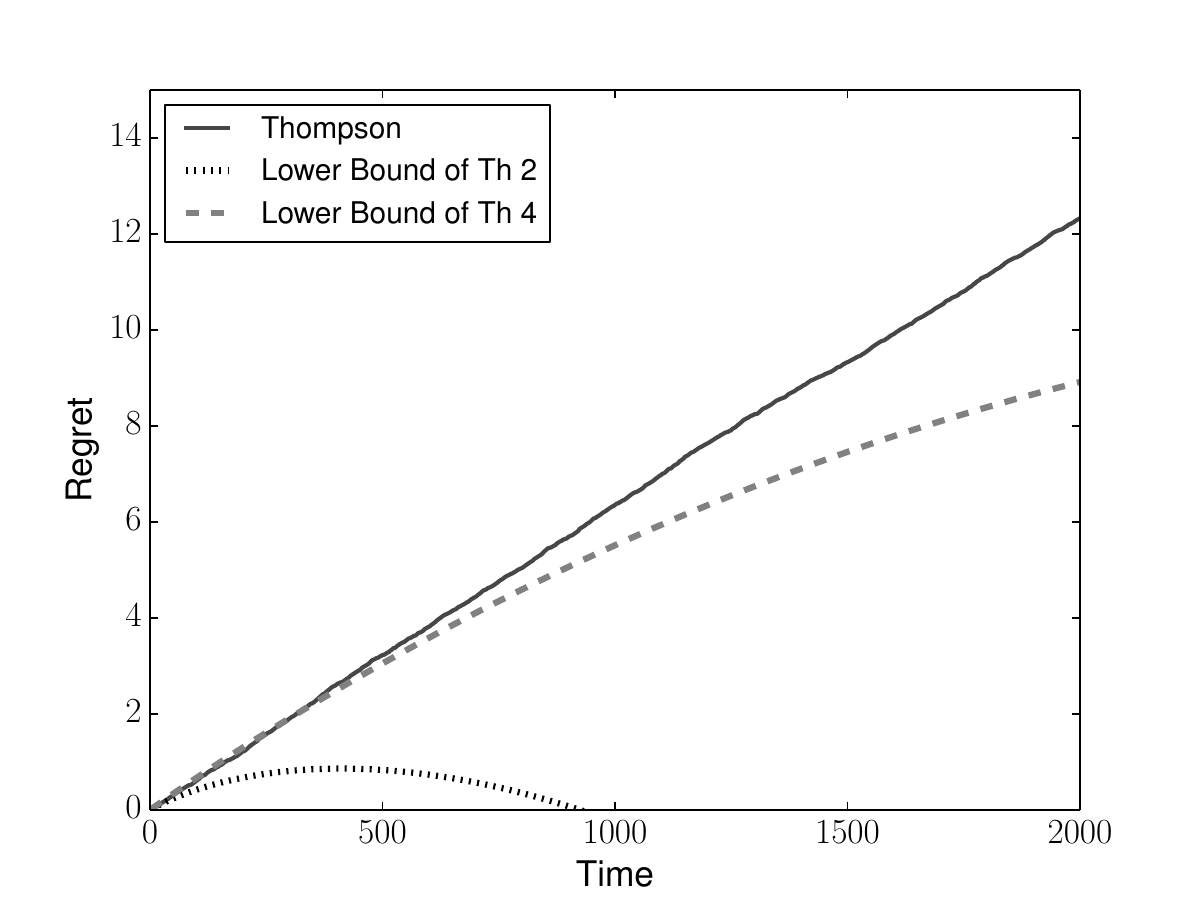} &
\includegraphics[width=0.47\textwidth]{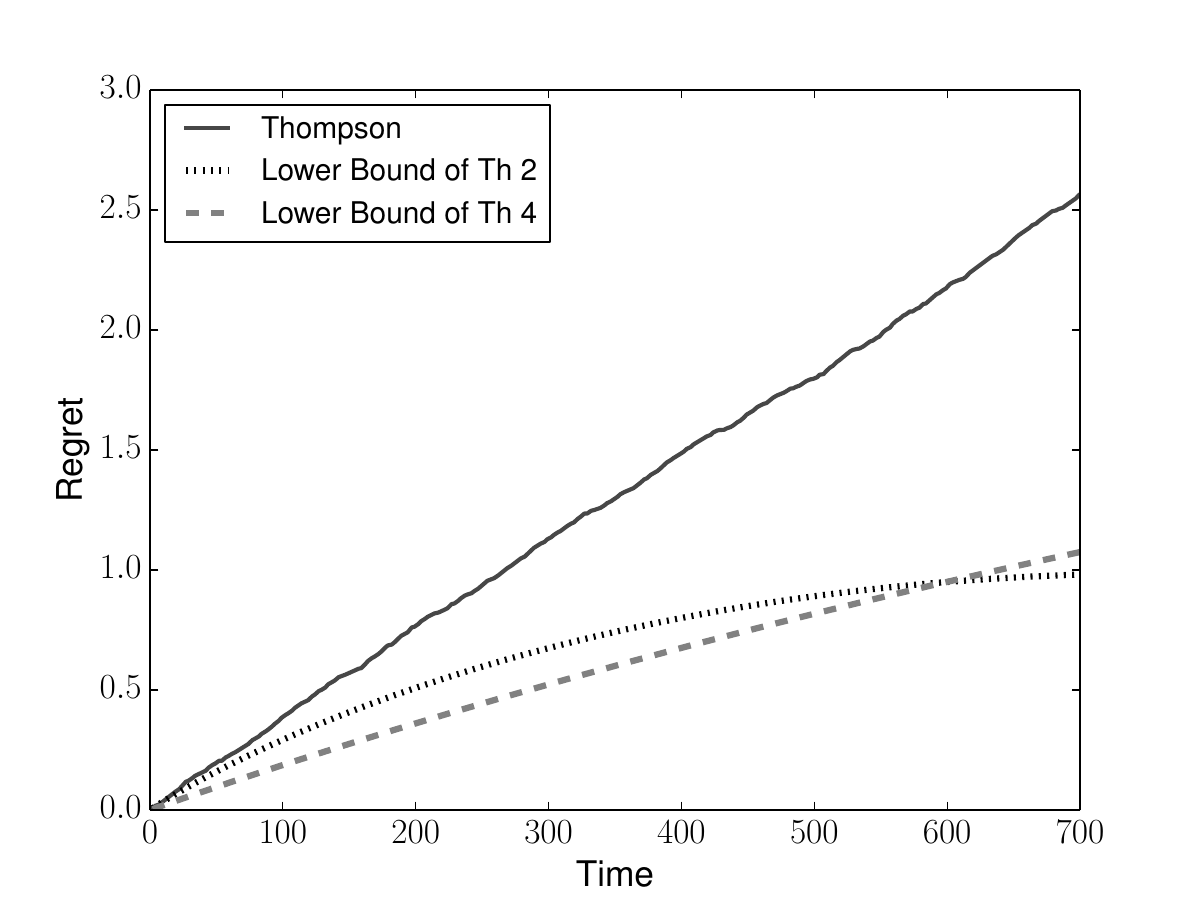}
\end{tabular}
\caption{\label{fig:th4}
Expected regret of \citemor{Th33} Sampling (\emph{blue, solid} line) on Bernoulli bandit problems,
versus the lower bound (\emph{red, dashed} line) of Theorem~\ref{th:4} using~\eqref{eq:LBth4} and the lower bound (\emph{black, dotted} line) of Theorem~\ref{th:2} using~\eqref{eq:towerrule},
for the model $\cD$ of all Bernoulli distributions;
expectations are approximated over $3,000$ runs.
\smallskip \newline
\emph{Left}: parameters $(\mu_a)_{1 \leq a \leq 10} =(0.05,\, 0.043, \,\ldots,\, 0.043)$,
with characteristic time $K/\big(8 \, A^\star_{\unu} \, \cK^{\max}_{\unu}\big) \approx 1,250$.
\newline
\emph{Right}: parameters $(\mu_a)_{1 \leq a \leq 7} =
(0.05, \, 0.048,\, 0.047, \, 0.046, \, 0.045, \, 0.044, \, 0.043)$, with
$K/\big(8 \, A^\star_{\unu} \, \cK^{\max}_{\unu}\big) \approx 1,619$.}
\end{figure}

Theorem~\ref{th:4} is illustrated in Figure~\ref{fig:th4}.
For a large number of arms, the regret lower bound~\eqref{eq:LBth4}
deriving as a consequence of the considered theorem is larger
than a bound based on the decomposition of the regret~\eqref{eq:towerrule}
and Theorem~\ref{th:2}.

\section{Non-asymptotic bounds for large T.}
\label{sec:largeT}

We restrict our attention to well-behaved models and uniformly super-fast convergent strategies.
For a given model $\cD$, we denote by $E(\cD)$ the interior of the set of all expectations
of distributions in $\cD$. That a model is well-behaved means that
the function $\Kinf$ is locally Lipschitz continuous in its second variable,
as is made formal in the following definition.

\begin{defi}
\label{eq:wbm}
A model $\mathcal{D}$ is well behaved if
there exist two functions $\varepsilon_{\cD} : E(\cD) \to (0,+\infty)$
and $\omega_{\cD} : \cD \times E(\cD) \to (0,+\infty)$
such that for all distributions $\nu_a \in \cD$ and
all $x \in E(\cD)$ with $x > E(\nu_a)$,
\[
\forall \eps < \eps_{\cD}(x), \qquad
\Kinf(\nu_a,x+\eps,\cD) \leq \Kinf(\nu_a,x,\cD) + \eps \, \omega_{\cD}(\nu_a,x)\,.
\]
\end{defi}

We could have considered a more general definition, where the upper bound would
have been any vanishing function of $\eps$, not only a linear function of $\eps$. However,
all examples considered in this paper (see Section~\ref{sec:wbm}) can be associated with such a linear difference.
Those examples of well-behaved models include parametric families like regular exponential families,
as well as more massive classes, like the set of all distributions with bounded support
(with or without a constraint on the finiteness of support).
Some of these examples, namely, regular exponential families and finitely-supported distributions with common
bounded support, were the models studied in~\citemor{klucb}
to get non-asymptotic upper bounds on the regret of the optimal order~\eqref{eq:obj}.

\begin{defi}
\label{eq:superconsistency}
A strategy $\psi$ is uniformly super-fast convergent on a model $\mathcal{D}$ if there exists a constant $C_{\psi,\mathcal{D}}$
such that for all bandit problems $\unu$ in $\mathcal{D}$, for all
suboptimal arms $a$, for all $T \geq 2$,
\[
\E_{\unu}\bigl[N_{\psi,a}(T)\bigr] \leq C_{\psi,\mathcal{D}} \frac{\ln T}{\Delta_a^2}\,.
\]
\end{defi}

Uniform super-fast convergence is a refinement of the notion of uniform fast convergence based on two considerations.
First, that there exist such strategies, for instance,
the \algo{UCB} strategy of \citemor{AuCBFi02} on any bounded model $\cD$, i.e.,
a model with distributions all supported within a common bounded interval $[m,M]$.
Second, Pinsker's inequality (see Appendix~\ref{sec:reminderinfotheory}) and Lemma~\ref{lm:DPE} entail in particular
that for such bounded models $\cD$,
\[
\Kinf(\nu_a,\mu^\star,\cD) \geq \kl\!\left( \frac{\mu_a-m}{M-m}, \, \frac{\mu^\star-m}{M-m}\right)
\geq \frac{2}{(M-m)^2} \Delta_a^2\,;
\]
therefore, the upper
bound stated in the definition of uniform super-fast convergence is still
weaker than the lower bound~\eqref{eq:obj}.

Note that Definition~\ref{eq:superconsistency} could be relaxed even more: we are mostly interested
therein in the logarithmic growth rate $\ln T$. We imposed the $C_{\psi,\cD}/\Delta_a^2$ upper bound mostly for simplicity
and readability of the calculations that lead to Theorem~\ref{th:LBlargeT-NA}.
It would be of course possible to rather consider more abstract problem-dependent constants of the
form $C_{\psi,\cD}(a,\nu)$, at least as soon as some minimal properties are assumed with respect to the behavior
of such constants as functions of the gap $\mu^\star - \mu_a$.

\subsection{A general non-asymptotic lower bound.}

Throughout this subsection, we fix a strategy $\psi$ that is uniformly super-fast convergent
with respect to a model $\mathcal{D}$.
We recall that we denote by $\cA^\star(\unu)$ the set of optimal arms of the bandit problem $\unu$
and let $A^\star_{\unu}$ be its cardinality.
We adapt the bounds~\eqref{eq:funda} and~\eqref{eq:BuKa-almost2} by using this time
\[
Z = \frac{1}{T} \sum_{a^\star \in \cA^\star(\unu)} N_{\psi,a^\star}(T)
\]
and $\kl(p,q) \geq p \ln(1/q) - \ln 2$, see~\eqref{ineg:kl}.
For all bandit problems $\unu'$
that only differ from $\unu$ as far a suboptimal arm $a$ is concerned, whose
distribution of payoffs $\nu'_a \in \mathcal{D}$ is such that $\mu'_a = E(\nu'_a) > \mu^\star$,
we get
\begin{equation}
\label{eq:nonasympt}
\E_{\unu} \! \bigl[ N_{\psi,a}(T) \bigr]
\geq \frac{1}{\KL(\nu_a,\nu'_a)} \left( \E_{\unu}[Z] \,
\ln \frac{1}{\E_{\unu'}[Z]} - \ln 2 \right).
\end{equation}
We restrict our attention to distributions $\nu'_a \in \mathcal{D}$
such that the gaps for $\unu'$ associated with optimal arms $a^\star \in \cA^\star(\unu)$
of $\unu$ satisfy $\uDelta = \mu'_a - \mu^\star \geq \eps$,
for some parameter $\eps > 0$ to be defined  by the analysis.
By uniform super-fast convergence, on the one hand,
\[
\E_{\unu}[Z] = 1 - \frac{1}{T} \sum_{a \not\in \cA^\star(\unu)} \E_{\unu}\bigl[N_{\psi,a}(T)\bigr]
\geq 1 - \frac{1}{T} \left( C_{\psi,\mathcal{D}}
\sum_{a \not\in \cA^\star(\unu)} \frac{1}{\Delta_a^2} \ln T \right);
\]
on the other hand,
\[
\E_{\unu'}[Z] = \frac{1}{T} \sum_{a^\star \in \cA^\star(\unu)} \E_{\unu'}\bigl[N_{\psi,a}(T)\bigr]
\leq \frac{A^\star_{\unu}\, C_{\psi,\mathcal{D}}}{\uDelta^2} \frac{\ln T}{T}\,.
\]
Denoting
\begin{equation}
\label{eq:Hnu}
H(\unu) = \sum_{a \not\in \cA^\star(\unu)} \frac{1}{\Delta_a^2}
\end{equation}
and using that $\uDelta \geq \eps$, a substitution of the two super-fast convergence inequalities
into~\eqref{eq:nonasympt} and an optimization over the considered distributions $\nu'_a$ leads to
\begin{equation}
\label{eq:nonasympt2opti}
\E_{\unu} \! \bigl[ N_{\psi,a}(T) \bigr]
\geq \frac{1}{\Kinf(\nu_a,\mu^\star+\eps,\cD)} \left(
1 - C_{\psi,\mathcal{D}} H(\unu) \frac{\ln T}{T} \right)
\ln \frac{T \eps^2}{\displaystyle{A^\star_{\unu}\, C_{\psi,\mathcal{D}} \ln T}}
- \frac{\ln 2}{\Kinf(\nu_a,\mu^\star+\eps,\cD)}\,.
\end{equation}
The obtained bound holds for all $T \geq 2$ (as in the definition of uniform super-fast convergence); however,
for small values of $T$, it might be negative, thus useless.

To proceed, we use the fact that the model $\mathcal{D}$ is well-behaved to relate
$\Kinf(\nu_a,\mu^\star+\eps,\cD)$ to $\Kinf(\nu_a,\mu^\star,\cD)$.
Since $1/(1+x) \geq 1-x$ for all $x \geq 0$, we get by Definition~\ref{eq:wbm}
\[
\forall \eps < \eps_{\cD}(\mu^\star), \qquad
\frac{1}{\Kinf(\nu_a,\mu^\star+\eps,\cD)} \geq \frac{1}{\Kinf(\nu_a,\mu^\star,\cD)}
\left( 1 - \eps \, \frac{\omega_{\cD}(\nu_a,\mu^\star)}{\Kinf(\nu_a,\mu^\star,\cD)} \right).
\]
Now, we set $\eps = \eps_T = (\ln T)^{-4}$.
Many other choices would have been possible, but this one is such that
$\eps_T \leq 0.0005$ already for $T \geq 1\,000$.
Putting all things together, from~\eqref{eq:nonasympt2opti},
from the fact that $(1-a)(1-b)(1-c) \geq 1 - (a+b+c)$ when
$0 \leq a,b,c \leq 1$, and from the bound $A^\star_{\unu} \leq K$,
we get the following theorem.

\begin{theorem}
\label{th:LBlargeT-NA}
For all uniformly super-fast convergent strategies $\psi$ on well-behaved models $\mathcal{D}$,
for all bandit problems $\unu$ in $\mathcal{D}$,
for all suboptimal arms $a$,
\begin{equation}
\label{eq:NAfinal}
\E_{\unu} \! \bigl[ N_{\psi,a}(T) \bigr]
\geq \frac{\ln T}{\Kinf(\nu_a,\mu^\star,\cD)} \bigl( 1
- (a_T + b_T + c_T) \bigr) - \frac{\ln 2}{\Kinf(\nu_a,\mu^\star,\cD)}\,,
\end{equation}
for all $T \geq 2$ large enough so that $(\ln T)^{-4} < \varepsilon_{\cD}(\mu^\star)$ and
\[
a_T = \frac{\omega_{\cD}(\nu_a,\mu^\star)}{\Kinf(\nu_a,\mu^\star,\cD)} (\ln T)^{-4}\,, \qquad
b_T = C_{\psi,\mathcal{D}} H(\unu) \frac{\ln T}{T}\,, \qquad
c_T = \frac{\ln \bigl( K\, C_{\psi,\mathcal{D}} (\ln T)^9\bigr)}{\ln T}\,,
\]
are all smaller than $1$, where $H(\unu)$ was defined in~\eqref{eq:Hnu}.
\end{theorem}

\begin{rema}
We have $(a_T + b_T + c_T) \ln T = \O\bigl(\ln(\ln T)\bigr)$.
The non-asymptotic bound~\eqref{eq:NAfinal} is therefore of the form
\[
\E_{\unu} \! \bigl[ N_{\psi,a}(T) \bigr]
\geq \frac{\ln T}{\Kinf(\nu_a,\mu^\star,\cD)} - \O\bigl(\ln(\ln T)\bigr)\,.
\]
Note that the second-order term of
typical non-asymptotic upper bounds (e.g., by \citemorpar{klucb})
had long been of the form $+ (\ln T)^\alpha$ for some $\alpha \in (0,1)$.
But recently, \citemoradd{HT15}{Theorem~5} showed that at least for models containing distributions
that have each a bounded support, the second-order is of order $- \ln(\ln T)$.
Our lower bound above thus shows the optimality of the order
of magnitude of this second-order term.
\end{rema}

\subsection{Two (and a half) examples of well-behaved models.}
\label{sec:wbm}

We consider first distributions with common bounded support (and the subclass of
such distributions with finite support); and then, regular exponential families.
The latter and the subclass of distributions with finite and bounded support
are the two models for which \citemor{klucb} could prove
non-asymptotic upper bounds matching the lower bound~\eqref{eq:obj}.

\paragraph{Distributions with common bounded support.}
We denote by $\mset$ the set of all probability distributions
over $[0,M]$, equipped with its Borel $\sigma$--algebra,
and restrict our model to such distributions with expectation
not equal to $M$.

\begin{lemma}
\label{lm:mset}
In the model $\mathcal{D} = \Bigl\{ \m \in \mset : E(\m) < M \Bigr\}$, we have
\begin{multline*}
\forall \m \in \mathcal{D}, \ \ \ \forall \mu^\star \in [0,M), \ \ \ \forall \eps \in \bigl( 0,(M-\mu^\star)/2 \bigr), \\
\qquad \qquad \Kinf(\m,\mu^\star+\eps,\cD) \leq \Kinf(\m,\mu^\star,\cD) - \ln \! \left( 1 -
\frac{2\eps}{M-\mu^\star} \right) .
\end{multline*}
In particular, for all $\m \in \mathcal{D}$ and $\mu^\star \in [0,M)$,
\[
\forall \eps \in \bigl( 0,(M-\mu^\star)/4 \bigr),
\quad \qquad \Kinf(\m,\mu^\star+\eps,\cD) \leq \Kinf(\m,\mu^\star,\cD) +
\frac{4\eps}{M-\mu^\star}\,.
\]
\end{lemma}

\begin{MORproof}
We fix $\m$, $\mu^\star$ and $\eps$ as indicated for the first bound; in particular, $\mu^\star+\eps < M$.
Since $\m$ is a probability distribution, it has at most countably many atoms; therefore,
there exists some $x \in (\mu^\star+\eps,M)$ such that $\m(\{x\}) = 0$
and $x \geq (M+\mu^\star)/2$. In particular,
$\m$ and the Dirac measure $\delta_x$ at this point are singular measures.

We consider some $\m' \in \mathcal{D}$ such that $E(\m') > \mu^\star$
and $\m \ll \m'$ (i.e., $\m$ is absolutely continuous with respect to $\m'$).
Such distributions exist and they are the only interesting ones in the defining
infimum of $\Kinf(\m,\mu^\star,\cD)$. We associate with $\m'$ the distribution
\[
\m'_\alpha = (1-\alpha) \m' + \alpha \delta_x\,,
\qquad \mbox{for the value} \qquad \alpha = \frac{\eps}{x-\mu^\star} \in (0,1)\,.
\]
The expectation of $\m'_\alpha$ satisfies
\begin{equation}
\label{eq:expmalpha}
E\bigl( \m'_\alpha \bigr) > (1-\alpha) \mu^\star + \alpha x = \mu^\star + \alpha (x-\mu^\star) = \mu^\star+\eps\,.
\end{equation}
Now, $\m \ll \m'$ entails that $\m \ll \m'_\alpha$ as well, with respective densities
satisfying (because $\m$ and $\delta_x$ are singular)
\[
\frac{\d\m}{\d\m'_\alpha} = \frac{1}{1-\alpha}\,\,\frac{\d\m}{\d\m'}
\qquad \mbox{and} \qquad
\frac{\d\m}{\d\m'_\alpha}(x) = 0\,.
\]
Therefore,
\[
\KL(\m,\m'_\alpha) = \bigintsss \left( \ln \frac{\d\m}{\d\m'_\alpha} \right)\d\m
= \ln \frac{1}{1-\alpha} + \bigintsss \left( \ln \frac{\d\m}{\d\m'} \right)\d\m
= \ln \frac{1}{1-\alpha} + \KL(\m,\m')\,.
\]
Since $\alpha$ decreases with $x$ and $x \geq (M+\mu^\star)/2$, we get
$\alpha \leq 2\eps/(M-\mu^\star)$. We substitute this bound
in the inequality above and take the infimum in both sides,
considering~\eqref{eq:expmalpha}, to get the first claimed bound.
The second bound follows from the inequality $-\ln(1-x) \leq 2x$ for $x \in [0,1/2]$.
\end{MORproof}

\begin{rema}
We denote by $\msetfin$ the subset of $\mset$ formed by probability distributions with
finite support.
The proof above shows that the bound of Lemma~\ref{lm:mset} also holds
for the model
\[
\mathcal{D} = \Bigl\{ \m \in \msetfin : E(\m) < M \Bigr\}\,.
\]
\end{rema}

\paragraph{Regular exponential families.}
Another example of well-behaved models is given by regular exponential families,
see \citemor{LeCa98} for a thorough exposition or \citemor{klucb} for an alternative
exposition focused on multi-armed bandit problems. \medskip

Such a family $\mathcal{D}$ is indexed by an open set $I = (m,M)$, where for each $\mu \in I$
there exists a unique distribution $\nu_\mu \in \mathcal{D}$ with expectation $\mu$.
(The bounds $m$ and $M$ can be equal to $\pm\infty$.) A key property of such a family is
that the Kullback-Leibler divergence between two of its elements can be represented\footnote{This function $g$
has an intrinsic definition as the convex conjugate
of the log-normalization function $b$ in the natural parameter space $\Theta$,
where $b$ can also be seen as a primitive of the expectation function $\Theta \to I$. But these
properties are unimportant here.} by
a twice differentiable and strictly convex function $g : I \to \R$, with increasing first
derivative $\dot{g}$ and continuous second derivative $\ddot{g} \geq 0$, in the sense that
\begin{equation}
\label{eq:KLexp}
\forall \, (\mu,\mu') \in I^2, \qquad \KL\bigl(\nu_\mu,\nu_{\mu'}\bigr)
= g(\mu) - g(\mu') - (\mu-\mu') \, \dot{g}(\mu')\,.
\end{equation}
In particular, $\mu' \mapsto \KL\bigl(\nu_\mu,\nu_{\mu'}\bigr)$
is strictly convex on $I$, thus is increasing on $[\mu,M)$. This entails
that
\begin{equation}
\label{eq:KinfKLexp}
\forall \, (\mu,\mu^\star) \in I^2 \ \ \mbox{s.t} \ \ \mu < \mu^\star, \qquad
\Kinf(\nu_\mu,\mu^\star,\cD) = \KL\bigl(\nu_\mu,\nu_{\mu^\star}\bigr)\,.
\end{equation}
In the lemma below, we restrict our attention to $\eps > 0$ such that
$\mu^\star+\eps \in I$, e.g., to $\eps < B_{\mu^\star}$ where
\begin{equation}
\label{eq:Bmustar}
B_{\mu^\star} = \min\!\left\{ \frac{M-\mu^\star}{2}, \, 1 \right\}.
\end{equation}
The minimum with $1$ is considered merely for $B_{\mu^\star}$ to always
have a finite value; otherwise, the bound in the lemma below would
be uninformative.

\begin{lemma}
In a model $\mathcal{D}$ given by a regular exponential family indexed by $I = (m,M)$
and whose Kullback-Leibler divergence~\eqref{eq:KLexp} is represented by a function $g$,
we have, with the notation~\eqref{eq:Bmustar},
\[
\forall \, \mu < \mu^\star \ \mbox{\rm of} \ I, \ \ \ \forall \, 0 < \eps
< B_{\mu^\star},
\quad \qquad \Kinf(\nu_\mu,\mu^\star+\eps,\cD) \leq \Kinf(\nu_\mu,\mu^\star,\cD)
+ \eps \, \bigl( \mu^\star + B_{\mu^\star} - \mu \bigr) \, G_{\mu^\star}
\]
where $G_{\mu^\star} = \max\bigl\{ \ddot{g}(x) : \ \ \mu^\star \leq x \leq \mu^\star + B_{\mu^\star} \bigr\}$.
\end{lemma}

\begin{MORproof}
Since $\mu < \mu^\star$, we get by~\eqref{eq:KLexp}
and~\eqref{eq:KinfKLexp}
\begin{eqnarray*}
\lefteqn{\Kinf(\nu_\mu,\mu^\star+\eps,\cD) - \Kinf(\nu_\mu,\mu^\star,\cD)} \\
& = & g(\mu^\star) - g(\mu^\star+\eps) - \bigl( \mu - (\mu^\star+\eps) \bigr) \, \dot{g}(\mu^\star+\eps)
+ (\mu-\mu^\star) \, \dot{g}(\mu^\star) \\
& = & \underbrace{g(\mu^\star) - g(\mu^\star+\eps) + \eps\,\dot{g}(\mu^\star)}_{\leq 0}
+ \bigl( (\mu^\star+\eps) - \mu \bigr) \bigl( \dot{g}(\mu^\star+\eps) - \dot{g}(\mu^\star) \bigr)\,,
\end{eqnarray*}
where the inequality is obtained by convexity of $g$. The proof is concluded by an application of the
mean-value theorem,
\[
\dot{g}(\mu^\star+\eps) - \dot{g}(\mu^\star) \leq \eps \, \max_{(\mu^\star,\mu^\star+\eps)} \ddot{g}\,,
\]
and the bound $\eps \leq B_{\mu^\star}$.
\end{MORproof}

The upper bound obtained on $\Kinf(\nu_\mu,\mu^\star+\eps,\cD) - \Kinf(\nu_\mu,\mu^\star,\cD)$
equals $\eps \, \bigl( \mu^\star + B_{\mu^\star} - \mu \bigr) \,
G_{\mu^\star}$. The examples below propose concrete upper bounds for $G_{\mu^\star}$
in different exponential families. None of these upper bounds actually involves $B_{\mu^\star}$
as various monotonicity arguments can be invoked.

\begin{example}
For Poisson distributions, we have $I = (0,+\infty)$ and
\[
\KL\bigl(\nu_\mu,\nu_{\mu'}\bigr) = \mu' - \mu + \mu \ln \frac{\mu}{\mu'}\,.
\]
We may take $g(\mu) = \mu \ln \mu - \mu$, so that $\ddot{g}(\mu) = 1/\mu$ and
$G_{\mu^\star} = 1/\mu^\star$.
\end{example}

\begin{example}
For Gamma distributions with known shape parameter $\alpha > 0$
(e.g., the exponential distributions when $\alpha = 1$), we have $I = (0,+\infty)$ and
\[
\KL\bigl(\nu_\mu,\nu_{\mu'}\bigr) = \alpha \left( \frac{\mu}{\mu'} - 1 - \ln \frac{\mu}{\mu'} \right).
\]
We may take $g(\mu) = - \alpha \ln \mu$, so that $\ddot{g}(\mu) = \alpha/\mu^2$ and
$G_{\mu^\star} = \alpha/(\mu^\star)^2$.
\end{example}

\begin{example}
For Gaussian distributions with known variance $\sigma^2 > 0$,
we have $I = (0,+\infty)$ and
\[
\KL\bigl(\nu_\mu,\nu_{\mu'}\bigr) = \frac{(\mu-\mu')^2}{2\sigma^2}\,.
\]
We may take $g(\mu) = \mu^2/(2\sigma^2)$, so that $\ddot{g}(\mu) = 1/\sigma^2$ and
$G_{\mu^\star} = 1/\sigma^2$.
\end{example}

\begin{example}
For binomial distributions for $n$ samples
(e.g., Bernoulli distributions when $n=1$),
we have $I = (0,n)$ and
\[
\KL\bigl(\nu_\mu,\nu_{\mu'}\bigr) = \mu \ln \frac{\mu}{\mu'} + (n-\mu) \ln \frac{n-\mu}{n-\mu'}\,.
\]
We may take $g(\mu) = \mu \ln \mu + (n-\mu)\ln(n-\mu)$,
so that $\ddot{g}(\mu) = n/\bigl( \mu(n-\mu) \bigr)$. A possible upper bound is
\[
G_{\mu^\star} \leq \frac{2n}{\mu^\star(n-\mu^\star)}\,.
\]
This can be seen by noting that $B_{\mu^\star} \leq (n-\mu^\star)/2$ so
that any $\mu \in [\mu^\star, \, \mu^\star+B_{\mu^\star}]$
is such that
$\mu \geq \mu^\star$ and
$n-\mu \geq n - \mu^\star - B_{\mu^\star} \geq (n-\mu^\star)/2$.
\end{example}

\begin{APPENDICES}
\section{Reminder of some elements of information theory.}
\label{sec:reminderinfotheory}

For the sake of self-completeness we recall two selected
basic facts pertaining to Kullback-Leibler divergences.

\paragraph{The data-processing inequality.} The most elegant proof we are aware of
relies on a conditional Jensen's inequality applied to $t \mapsto t \ln t$;
see~\citemor{DPineq}.

\begin{lemma}
\label{lm:DPineq}
Consider a measurable space $(\Gamma,\mathcal{G})$ equipped with two distributions
$\P_1$ and $\P_2$, any other $(\Gamma',\mathcal{G}')$ measurable space,
and any random variable $X : (\Gamma,\mathcal{G}) \to (\Gamma',\mathcal{G}')$.
Then,
\[
\KL\bigl(\P_1^X,\P_2^X\bigr) \leq \KL(\P_1,\P_2)\,,
\]
where $\P_1^X$ and $\P_2^X$ denote the respective distributions
of $X$ under $\P_1$ and~$\P_2$.
\end{lemma}

\paragraph{On local refinements of Pinsker's inequality.}
Pinsker's inequality reads, for Bernoulli distributions, in its most classical form:
\begin{equation}
\label{eq:Pinskerclassical}
\forall (p,q) \in [0,1]^2, \qquad \kl(p,q) \geq 2(p-q)^2\,.
\end{equation}
The lemma below offers a local refinement of Pinsker's inequality for Bernoulli distributions;
the classical form~\eqref{eq:Pinskerclassical} follows by noting that $x(1-x) \leq 1/4$ for $x \in [0,1]$.
\citemoradd{klucb}{Lemma~3 in Appendix A.2.1} offer an extension of this local refinement
to any one-parameter regular exponential family.
\begin{lemma}
\label{lm:Pinskerlocal}
For $0 \leq p < q \leq 1$, we have \quad
$\displaystyle{
\kl(p,q) \geq \frac{1}{2 \displaystyle{\max_{x \in [p,q]}} x(1-x)}(p-q)^2
\geq \frac{1}{2 q}(p-q)^2\,.
}$
\end{lemma}

\begin{MORproof}
We may assume that $p > 0$ and $q < 1$, since for $p=0$, the result follows by continuity,
and for $q =1$, the inequality is void, as $\kl(p,1) = +\infty$ when $p < 1$.
The first and second derivative of $\kl$ equal
\[
\frac{\partial}{\partial p}\kl(p,q) = \ln p - \ln(1-p) - \ln q  + \ln(1-q)
\quad \mbox{and} \quad
\frac{\partial^2}{\partial^2 p}\kl(p,q) = \frac{1}{p} + \frac{1}{1-p} = \frac{1}{p(1-p)}\,.
\]
By Taylor's equality, there exists $r \in [p,q]$ such that
\[
\kl(p,q) = \underbrace{\kl(q,q)}_{= 0} + (p-q) \underbrace{\frac{\partial}{\partial p}\kl(q,q)}_{=0}
+ \frac{(p-q)^2}{2} \underbrace{\frac{\partial^2}{\partial^2 p}\kl(r,q)}_{=1/(r(1-r))}\,.
\]
The proof of the first inequality is concluded by upper bounding $r(1-r)$ by $\displaystyle{\max_{x \in [p,q]} x(1-x)}$. \\
The second inequality follows from
$\displaystyle{\max_{x \in [p,q]} x(1-x) \leq \max_{x \in [p,q]} x \leq q}$.
\end{MORproof}

\section{Re-derivation of other earlier lower bounds}
\label{sec:appLB}

In this section, we re-derive the bounds discussed in Section~\ref{sec:other-blb},
based on our fundamental inequality~\eqref{eq:funda}. We do so to illustrate the power and the versatility
of~\eqref{eq:funda}. However, we point out again that the lower bounds discussed here are
much weaker than the ones derived in the main body of the paper: in the terminology
of Section~\ref{sec:other-blb}, they are of the form (well-chosen) rather than of the form~(all).

\subsection{Distribution-free lower bound.}
\label{sec:rederiv-distrfree}

We consider the bound~\eqref{eq:distrfree} recalled in Section~\ref{sec:other-blb}.
More specifically, we re-prove Theorem~A.2 of \citemor{AuCBFrSc02},
from which the stated bound~\eqref{eq:distrfree} follows by optimization
over $\varepsilon$.

\begin{theorem}
\label{th:distr-free}
Consider the bandit model $\mathcal{D} = \mathcal{M}\bigl([0,1]\bigr)$
of all probability distributions over $[0,1]$.
For all $\varepsilon \in (0,1/2)$,
for all strategies $\psi$,
there exists a bandit problem $\unu'$ in $\mathcal{M}\bigl([0,1]\bigr)$ such that
\[
R_{\psi,\unu',T} \geq T \varepsilon \left( 1 - \frac{1}{K} - \frac{1}{2} \sqrt{\frac{T}{K}
\ln \frac{1}{1-4\varepsilon^2}} \right).
\]
This problem $\unu'$ can be given by Bernoulli distributions, with parameters $1/2$ for
all arms but one, for which the parameter is $1/2+\varepsilon$.
\end{theorem}

As a consequence, the worst-case regret of any strategy $\psi$ against
all bandit problems $\unu$ in $\mathcal{M}\bigl([0,1]\bigr)$ is lower bounded as
announced in~\eqref{eq:distrfree}:
\[
\sup_{\unu} R_{\psi,\unu,T} \geq \sup_{\varepsilon \in (0,1/2)}
T \varepsilon \left( 1 - \frac{1}{K} - \frac{1}{2} \sqrt{\frac{T}{K}
\ln \frac{1}{1-4\varepsilon^2}} \right) \geq
\frac{1}{20} \min\Bigl\{ \sqrt{KT}, \, T \Bigr\}\,.
\]
The second inequality above is proved by a simple calculation indicated
after the proof of Theorem~A.2 of~\citemor{AuCBFrSc02}: pick $\varepsilon = \min\bigl\{\sqrt{K/T},\,1\bigr\}/4$
and use $-\ln(1-u) \leq \bigl(4\ln(4/3)\bigr)u$ for $u \in (0,1/4)$.
The constant $1/20$ can actually be improved into $1/8$,
see~\citemoradd{CBLu06}{Theorem~6.11}.

\begin{MORproof}
We fix a strategy and $\varepsilon \in (0,1/2)$.
We denote by $\unu$ the bandit problem where all distributions are given
by Bernoulli distributions with parameter $1/2$. There exists an arm $k \in \{1,\ldots,K\}$
such that $\E_{\unu}\!\bigl[N_{\psi,k}(T)\bigr] \leq T/K$, as these $K$ numbers of pulls sum up to $T$.
We define the bandit problem $\unu'$ by $\nu'_a = \nu_a$ for $a \ne k$, that is, $\nu'_a$
is a symmetric Bernoulli distribution, while $\nu'_k$ is the Bernoulli distribution with
parameter $1/2+\varepsilon$. By~\eqref{eq:towerrule}, we have
\begin{equation}
\label{eq:dfree1}
R_{\psi,\unu',T} = \sum_{a \ne k} \varepsilon \, \E_{\unu'}\bigl[N_{\psi,a}(T)\bigr]
= T \varepsilon \left( 1 - \frac{\E_{\unu'}\!\bigl[N_{\psi,k}(T)\bigr]}{T} \right).
\end{equation}
A direct computation of $\kl(1/2,\,1/2+\varepsilon)$ and the application of~\eqref{eq:funda}
indicate that
\[
\frac{\E_{\unu} \! \bigl[ N_{\psi,k}(T) \bigr]}{2} \ln \frac{1}{1-4\varepsilon^2}
= \E_{\unu} \! \bigl[ N_{\psi,k}(T) \bigr] \, \kl\biggl(\frac{1}{2},\,\frac{1}{2}+\varepsilon\bigg)
\geq \kl\Biggl( \frac{\E_{\unu} \! \bigl[ N_{\psi,k}(T) \bigr]}{T}, \,\, \frac{\E_{\unu'} \! \bigl[ N_{\psi,k}(T) \bigr]}{T} \Biggr)\,.
\]
Now, Pinsker's inequality (in its classical form, see Appendix~\ref{sec:reminderinfotheory}) ensures that
\[
\frac{\E_{\unu} \! \bigl[ N_{\psi,k}(T) \bigr]}{2} \ln \frac{1}{1-4\varepsilon^2}
\geq
\kl\Biggl( \frac{\E_{\unu} \! \bigl[ N_{\psi,k}(T) \bigr]}{T}, \,\, \frac{\E_{\unu'} \! \bigl[ N_{\psi,k}(T) \bigr]}{T} \Biggr)
\geq 2 \left( \frac{\E_{\unu'} \! \bigl[ N_{\psi,k}(T) \bigr]}{T} - \frac{\E_{\unu} \! \bigl[ N_{\psi,k}(T) \bigr]}{T}
\right)^{\!\! 2}.
\]
Solving for $\E_{\unu'} \! \bigl[ N_{\psi,k}(T) \bigr]/T$, based
on whether $\E_{\unu'} \! \bigl[ N_{\psi,k}(T) \bigr]/T$ is larger or smaller than
$\E_{\unu} \! \bigl[ N_{\psi,k}(T) \bigr]/T$, we get, in all cases,
\[
\frac{\E_{\unu'} \! \bigl[ N_{\psi,k}(T) \bigr]}{T} \leq
\frac{\E_{\unu} \! \bigl[ N_{\psi,k}(T) \bigr]}{T} +
\frac{1}{2} \sqrt{\E_{\unu} \! \bigl[ N_{\psi,k}(T) \bigr]
\ln \frac{1}{1-4\varepsilon^2}}\,.
\]
The proof is concluded by substituting the fact that
$\E_{\unu}\!\bigl[N_{\psi,k}(T)\bigr] \leq T/K$ by definition of $k$,
and by combining the obtained inequality with~\eqref{eq:dfree1}.
\end{MORproof}

The short proof above actually re-uses absolutely all the original arguments
of~\citemor{AuCBFrSc02}: the same Bernoulli distributions, the chain rule
for Kullback-Leibler divergences, Pinsker's inequality. It is merely stated
in a compact way, that puts under the same umbrella the distribution-dependent
and the distribution-free lower bounds for multi-armed bandit problems.

\subsection{Lower bounds for the case when $\mu^\star$ or the gaps~$\Delta$ are known.}
\label{sec:rederiv-BPR}

We consider here the second framework discussed in Section~\ref{sec:other-blb},
with sub-Gaussian bandit problems.
For simplicity, and following \citemor{BPR},
we restrict our attention to lower bounds
for two-armed bandit problems (i.e., for $K=2$).

\paragraph{Known largest expected payoff $\mu^\star$ but unknown gap $\Delta$.}

The lower bound stated in Theorem~\ref{th:LB-BPR} below corresponds to Theorem~8 of~\citemor{BPR},
later revisited by the authors, see~\cite{err}. It turns out that, as hinted at
in, e.g., \citemoradd{ESAIM}{end of Section~1.4},
the initially claimed $\ln T$ dependency is incorrect and a bounded
regret can be guaranteed.
As shown in Theorem~\ref{th:finiteRegret:mus:known} in the next section,
this bound on the regret can be as small as $\ln(1/\Delta)/\Delta$.
The lower bound we could get using our techniques is of order $1/\Delta$.

To state it, we restrict our attention to strategies $\psi$ symmetric in some sense,
e.g., in the sense of Definition~\ref{def:ps} stated later on.
We actually need very little symmetry here: the considered strategies
$\psi$ should just be such that in the bandit problem
$\unu_0 = \bigl( \cN(0,1), \, \cN(0,1) \bigr)$, in which the two arms
have the same distribution,
\begin{equation}
\label{eq:sym}
\E_{\unu_0} \! \bigl[ N_{\psi,1}(T) \bigr]
= \E_{\unu_0} \! \bigl[ N_{\psi,2}(T) \bigr] = \frac{T}{2}\,.
\end{equation}
Of course, all reasonable strategies are usually even more symmetric than that: they
are usually stable by permutations over the arms (i.e., they base their decisions only on the
payoffs received, not on the labeling of the arms).

\begin{theorem}
\label{th:LB-BPR}
For all $\Delta > 0$ we consider $\unu_\Delta = \bigl( \cN(0,1), \, \cN(-\Delta,1) \bigr)$
and $\unu_0 = \bigl( \cN(0,1), \, \cN(0,1) \bigr)$.
For all strategies $\psi$ that are symmetric in the sense of~\eqref{eq:sym},
for all $\Delta > 0$, for all $T \geq 1$,
\[
\E_{\unu_{\!\Delta}} \! \bigl[ N_{\psi,2}(T) \bigr] \geq \frac{1}{\Delta^2 + 1/T}
\qquad \mbox{and} \qquad
R_{\psi,\unu_{\!\Delta},T} \geq \frac{\Delta}{\Delta^2 + 1/T}\,.
\]
In addition, for all strategies $\psi$ and for all $T$ such that
$\E_{\unu_{\!\Delta}} \! \bigl[ N_{\psi,2}(T) \bigr] \geq 1$,
\[
\E_{\unu_{\!\Delta}} \! \bigl[ N_{\psi,2}(T) \bigr] \geq \min\left\{ \frac{2 \ln 2}{\Delta^2 + 2 \ln(4T)/T}, \,\, \frac{T}{2} \right\}
\qquad \mbox{and} \qquad
R_{\psi,\unu_{\!\Delta},T} \geq
\min\left\{ \frac{2 (\ln 2) \Delta}{\Delta^2 + 2 \ln(4T)/T}, \,\, \frac{T\Delta}{2} \right\}.
\]
\end{theorem}

Note that the constraint that $\E_{\unu_{\!\Delta}} \! \bigl[ N_{\psi,2}(T) \bigr] \geq 1$ is satisfied
for all $T \geq K$ by most of the reasonable strategies, as the latter typically start by playing each arm
once (in a random order).

\begin{MORproof}
We first note that $R_{\psi,\unu_{\!\Delta},T} = \Delta \,\, \E_{\unu_{\!\Delta}} \! \bigl[ N_{\psi,2}(T) \bigr]$.
Inequality~\eqref{eq:funda} entails that
\begin{multline}
\label{eq:LB-BPR-eq1}
\frac{\Delta^2}{2} \, \E_{\unu_{\!\Delta}}\!\bigl[N_{\psi,2}(T)\bigr] =
\E_{\unu_{\!\Delta}}\!\bigl[N_{\psi,2}(T)\bigr]\,\,\KL\bigl(\cN(-\Delta,1),\,\cN(0,1)\bigr) \\
\geq \kl\Biggl( \frac{\E_{\unu_{\!\Delta}} \! \bigl[ N_{\psi,2}(T) \bigr]}{T}, \,\, \frac{\E_{\unu_0} \! \bigl[ N_{\psi,2}(T) \bigr]}{T} \Biggr) =
\kl\Biggl( \frac{\E_{\unu_{\!\Delta}} \! \bigl[ N_{\psi,2}(T) \bigr]}{T}, \,\, \frac{1}{2} \Biggr),
\end{multline}
where we used respectively, for the two equalities, the closed-form expression
for the Kullback-Leibler divergences between Gaussian distribution with
the same variance and the symmetry assumption on the strategy.
Pinsker's inequality (in its classical form, see Appendix~\ref{sec:reminderinfotheory}),
followed by the inequality
\[
\forall \, x \in \R, \qquad 2 \left( \frac{1}{2} - x \right)^{\!\! 2}
\geq \frac{1}{2} - 2x\,,
\]
yields
\[
\frac{\Delta^2}{2} \, \E_{\unu_{\!\Delta}}\!\bigl[N_{\psi,2}(T)\bigr] \geq
2 \left( \frac{1}{2} - \frac{\E_{\unu_{\!\Delta}} \! \bigl[ N_{\psi,2}(T) \bigr]}{T} \right)^{\!\! 2}
\geq \frac{1}{2} - 2 \frac{\E_{\unu_{\!\Delta}} \! \bigl[ N_{\psi,2}(T) \bigr]}{T}\,.
\]
Simple manipulations entail the first claimed bound on $\E_{\unu_{\!\Delta}} \! \bigl[ N_{\psi,2}(T) \bigr]$.

For the second one, given the form of the lower bound, which involves a minimum with $T/2$,
it suffices to consider the case when $\E_{\unu_{\!\Delta}} \! \bigl[ N_{\psi,2}(T) \bigr]\,/T \leq 1/2$.
We use that
\[
\kl(x,\,1/2) = \ln 2 - h(x)\,, \qquad \mbox{where} \qquad h(x) = - \bigl( x \ln x + (1-x) \ln(1-x) \bigr)
\]
is the binary entropy function. Now, \citemoradd{PhDh}{page 8} indicates that $h(x) \leq x \ln(4/x)$
for all $x \in [0,1/2]$, so that, restricting our attention to $x \geq 1/T$, we get
\[
\forall \, x \in [1/T,\,1/2], \qquad \kl\bigg(\!x,\,\frac{1}{2}\bigg) \geq \ln 2 - x \ln\!\bigg(\frac{4}{x}\bigg) \geq \ln 2 - x \ln(4T)\,.
\]
Substituting this inequality into~\eqref{eq:LB-BPR-eq1}, using that
$x = \E_{\unu_{\!\Delta}} \! \bigl[ N_{\psi,2}(T) \bigr] \, /T$ lies in $[1/T,\,1/2]$,
concludes the proof.
\end{MORproof}

The proof above, which is simple and direct, illustrates the interest of Inequality~\eqref{eq:funda} over the standard approaches used so far to prove lower bounds in the same or similar settings.

\paragraph{Known gap $\Delta$ but unknown largest expected payoff $\mu^\star$.}

The lower bound stated in Theorem~\ref{th:LB-BPR2} below corresponds to Theorem~6 of~\citemor{BPR}.
It shows the optimality of the performance bound $\ln(T\Delta^2)/\Delta$
on the regret of the Improved--UCB strategy introduced by~\cite{ImprUCB} and further
studied by~\cite{ETC}. The latter improved the constant in the leading term,
which equals $\ln(T\Delta^2)/(2\Delta)$ when the gap $\Delta$ between the expected
payoffs between the two Gaussian arms with variance~$1$ is known.

We denote by $W$ the Lambert function: for all $u \geq 0$,
there exists a unique $v \geq 0$ such that $u\,\exp(u) = v$,
which is denoted by $v = W(u)$. The Lambert function $W$
is increasing on $[0,+\infty)$.
One may easily check that 
\[
\forall x \geq \mathrm{e}, \qquad \ln(x) - \ln\bigl(\ln(x)\bigr) \leq W(x) \leq \ln(x)\,.
\]

We state below two lower bounds: one for all strategies $\psi$, in terms of a maximum
between two regrets; and one for strategies that are symmetric and invariant by translation.
These properties of symmetry and invariance by translation are most natural requirements.
To define them, for all $c\in\R$ and all distributions $\nu$,
we denote by $\tau_c(\nu)$ the distribution of $Y+c$ when $Y\sim \nu$.

\begin{defi}
\label{def:invtrans:sym}
A strategy $\psi$ for $K$--armed bandits is symmetric and invariant by translation of the payoffs if for all permutations $\sigma$ of $\{1,\dots,K\}$, all $c \in \R$, and all $T\geq 1$, the distribution of $\big(N_{\psi,1}(T), \dots, N_{\psi,K}(T)\big)$ in the bandit problem $(\nu_1,\dots,\nu_K)$ is equal to the one of $\big(N_{\psi,\sigma^{-1}(1)}(T), \dots, N_{\psi,\sigma^{-1}(K)}(T)\big)$ in the bandit problem $\big(\tau_c(\nu_{\sigma(1)}),\dots,\tau_c(\nu_{\sigma(K)})\big)$.
\end{defi}

\begin{theorem}
\label{th:LB-BPR2}
We fix $\Delta > 0$ and consider $\unu_1 = \bigl( \cN(0,1), \, \cN(-\Delta,1) \bigr)$
and $\unu_2 = \bigl( \cN(0,1), \, \cN(\Delta,1) \bigr)$.
Then, for all strategies $\psi$, for all $T \geq 1$,
\begin{equation}
\label{eq:BPR2-us}
\max \bigl\{ R_{\psi,\unu_1,T}, \, R_{\psi,\unu_2,T} \bigr\} \geq \min \! \left\{
\frac{W\bigl(T \Delta^2/1.2\bigr)}{2 \Delta}, \,\, \frac{T\Delta}{2} \right\}.
\end{equation}
Or, alternatively, for all strategies $\psi$ that are symmetric and invariant by translation of the payoffs,
for all $T \geq 1$,
\[
R_{\psi,\unu_1,T} = R_{\psi,\unu_2,T} \geq \frac{W\bigl(T \Delta^2/1.2\bigr)}{2 \Delta}\,.
\]
\end{theorem}

\begin{remark}
We compare the obtained bound~\eqref{eq:BPR2-us} to Theorem~6 of~\citemor{BPR}.
First, the proof reveals that~\eqref{eq:BPR2-us} holds for all distributions
$\unu_1 = \bigl( P_0, \, \cN(-\Delta,1) \bigr)$
and $\unu_2 = \bigl( P_0, \, \cN(\Delta,1) \bigr)$
where $P_0$ is a probability distribution with expectation~$0$. For instance, \citemor{BPR}
considered the Dirac mass $\delta_0$ at~$0$.

Second, Theorem~6 of~\citemor{BPR} offers the bound
\begin{equation}
\label{eq:BPR2}
\max \bigl\{ R_{\psi,\unu_1,T}, \, R_{\psi,\unu_2,T} \bigr\} \geq \frac{\ln\bigl(T \Delta^2/2\bigr)}{4 \Delta}\,.
\end{equation}
Asymptotically, as $T \to +\infty$, our bound~\eqref{eq:BPR2-us} is smaller by a factor of $2$.
For small values of $T$ (or small values of $\Delta$), the bound~\eqref{eq:BPR2} is void
as the logarithmic term is non-positive, while our bound is always nonnegative.
The second argument of the minimum in~\eqref{eq:BPR2-us} is unimportant, as
the regret is always bounded by $T\Delta$.
\end{remark}

\begin{MORproof}
We have $R_{\psi,\unu_1,T} = \Delta \, \E_{\unu_{1}}\!\bigl[N_{\psi,2}(T)\bigr]$
and $R_{\psi,\unu_2,T} = \Delta \, \E_{\unu_{2}}\!\bigl[N_{\psi,1}(T)\bigr]$,
so that it suffices to lower bound
\[
x = \frac{1}{T} \max \Bigl\{ \E_{\unu_{1}}\!\bigl[N_{\psi,2}(T)\bigr], \,\,
\E_{\unu_{2}}\!\bigl[N_{\psi,1}(T)\bigr] \Bigr\}\,.
\]
We assume below that the maximum is
given by the first term; otherwise, the proof below should be adapted
by exchanging the roles of $\unu_1$ and $\unu_2$. Inequality~\eqref{eq:funda}
indicates that
\begin{multline*}
2 \, T \Delta^2 \, x =
2 \, \Delta^2 \, \E_{\unu_{1}}\!\bigl[N_{\psi,2}(T)\bigr] =
\E_{\unu_{1}}\!\bigl[N_{\psi,2}(T)\bigr]\,\,\KL\bigl(\cN(-\Delta,1),\,\cN(\Delta,1)\bigr) \\
\geq \kl\Biggl( \frac{\E_{\unu_{1}} \! \bigl[ N_{\psi,2}(T) \bigr]}{T}, \,\, \frac{\E_{\unu_2} \! \bigl[ N_{\psi,2}(T) \bigr]}{T} \Biggr) =
\kl\Biggl( \! x, \,\, 1 - \frac{\E_{\unu_2} \! \bigl[ N_{\psi,1}(T) \bigr]}{T} \Biggr)\,.
\end{multline*}
Given the form of the lower bound in the theorem, which involves a minimum
with $T\Delta/2$, we may assume, with no loss of generality, that $x \leq 1/2$.
Since $\kl(x,\,\cdot\,)$ is increasing on $[x,1]$ and since
\[
1 - \frac{\E_{\unu_2} \! \bigl[ N_{\psi,1}(T) \bigr]}{T}
\geq 1 - x \geq \frac{1}{2} \geq x\,,
\]
by definition of $x$ and the assumption $x \leq 1/2$, we get
\[
2 \, T \Delta^2 \, x \geq \kl(x,\,1-x) = (1-2x) \ln \frac{1-x}{x}\,.
\]
Note that the case $x = 0$ is excluded by the inequality above.
A function study shows that
\[
\forall \, x \in (0,1), \qquad (1-2x) \ln \frac{1-x}{x} \geq \ln \frac{1}{2.4\,x}\,.
\]
Substituting this lower bound and taking exponents, we are left with studying the inequality
\[
\exp \bigl( 2 \, T \Delta^2 \, x \bigr) \geq \frac{1}{2.4\,x}\,,
\qquad \mbox{or equivalently,} \qquad
2 \, T \Delta^2 \, x \,\, \exp \bigl( 2 \, T \Delta^2 \, x \bigr) \geq \frac{T \Delta^2}{1.2}\,.
\]
By definition of the Lambert function $W$, we rewrite this inequality as $2 \, T \Delta^2 \, x \geq
W\bigl(T \Delta^2/1.2\bigr)$,
which concludes the proof of the first statement.

For the second statement, we note that the property of invariance by translation of the payoffs
ensures that
\[
x = \frac{\E_{\unu_{1}}\!\bigl[N_{\psi,2}(T)\bigr]}{T} = \frac{\E_{\unu_{2}}\!\bigl[N_{\psi,1}(T)\bigr]}{T}\,.
\]
Therefore, the fundamental inequality~\eqref{eq:funda} directly gives in this case
\[
2 \, T \Delta^2 \, x \geq
\kl\Biggl(\frac{\E_{\unu_{1}} \! \bigl[ N_{\psi,2}(T) \bigr]}{T}, \,\, \frac{\E_{\unu_2} \! \bigl[ N_{\psi,2}(T) \bigr]}{T} \Biggr)
= \kl(x,\,1-x)\,,
\]
and we do not need to distinguish whether $x$ is larger than $1/2$ or not. The end of the proof
of the first statement of the theorem did not use that $x \leq 1/2$ and can still safely be followed
for the second statement.
\end{MORproof}

\section{A finite-regret algorithm when $\mu^\star$ is known.}
\label{sec:BPR}

In this section, and in this section only,
as we are discussing a specific strategy (described below in a box), we will not
index the regret, the number of times a given arm is pulled,
etc., by the said specific strategy.

We consider the sub-Gaussian framework described in Section~\ref{sec:other-blb}
and restrict our attention to the case when $\mu^\star$ is known. We provide a
refinement of the results of \citemoradd{BPR}{Section~3}, already known
by these authors themselves (see, e.g., \citemorpar{ESAIM}). The algorithm considered below is
inspired by Algorithm~1 of \citemor{BPR}. For each $t \geq 1$ and $a \in \{1,\ldots,K\}$ such
that $N_a(t) \geq 1$, we denote by
\[
\widehat{\mu}_{a,t} = \frac{1}{N_a(t)} \sum_{s=1}^t Y_s \, \ind_{ \{ A_s = a \}}
\]
the empirical mean of the rewards obtained between rounds $1$ and $t$
when playing arm $a$. \medskip

\begin{center}
\begin{algorithm}[H]
\label{algo_optimal_arm}
\caption{An algorithm with bounded regret, thanks to the knowledge of $\mu^\star$}
\SetKwInput{problem}{Bandit problem}
\SetKwInput{para}{Parameters}
\SetKwInput{for}{For}
\problem{$\unu = (\nu_a)_{a = 1,\ldots,K}$ where each $\nu_a$ is sub-Gaussian
in the sense of~\eqref{eq:sub-Gaussian}} \smallskip
\para{the value of $\displaystyle{\mu^\star = \max_{a = 1,\ldots,K} \mu_a}$} \smallskip
\for{each $t\in\{1,\ldots,K\}$, \textbf{do:} play arm $t$.} \smallskip
\for{each round $t \geq K+1$,}
\begin{enumerate}
\item Let $\displaystyle{\mathcal{C}_t = \left\{a \in \{1,\ldots,K\}\ : \ \ \ \widehat{\mu}_{a,t-1} - \mu^\star >
- \sqrt{\frac{4\ln N_a(t-1)}{N_a(t-1)}}\, \right\}}$ be the set of candidate arms; \smallskip
\item If $\mathcal{C}_t\neq\emptyset$, play an arm $A_t$ at random in $\mathcal{C}_t$, update $t:=t+1$; \smallskip
\item If $\mathcal{C}_t = \emptyset$, play $A_t=1, \,\, A_{t+1}=2, \,\, \ldots, \,\, A_{t+K} = t+K-1 $, update $t:=t+K$.
\end{enumerate}
\end{algorithm}
\end{center}
\smallskip

We use the notation introduced before~\eqref{eq:towerrule}, but, as indicated above, without the indexations in the
considered strategy.

\begin{theorem}\label{th:finiteRegret:mus:known}
For all bandit problems $\unu = (\nu_a)_{a = 1,\ldots,K}$ where each distribution $\nu_a$ is sub-Gaussian
in the sense of~\eqref{eq:sub-Gaussian}, the regret of the algorithm above is bounded by
\[
R_{\unu,T} \leq \sum_{a : \Delta_a > 0} \Biggl( \frac{36\ln(17/\Delta_a)}{\Delta_a} + 3\Delta_a \Biggr).
\]
\end{theorem}

\begin{MORproof}
We fix an optimal arm $a^\star$.
In view of~\eqref{eq:towerrule}, it suffices to bound $\E_{\unu} \! \bigl[ N_a(T) \bigr]$
for each suboptimal arm $a$. Each arm is played once between $1$ and $K$.
For all $t \geq K+1$, a suboptimal arm $a$ can only be played if $a \in \mathcal{C}_t$
(step~2 of the second for loop)
or if we are in a sequence where each arm is played successfully
(step~3 of the second for loop). In the latter case,
the set of candidate arms at round $t-a+1$ was empty. It did not contain $a^\star$.
This optimal arm is played also once in the sequence of pulls corresponding to step~3,
at time $t-a+a^\star+1$. At time $t-a+a^\star$ we still had
$N_{a^\star}(t-a+a^\star) = N_{a^\star}(t-a+1)$,
so that the condition for being a candidate was violated as well:
\[
\widehat{\mu}_{a^\star,t-a+a^\star} - \mu^\star \leq - \sqrt{\frac{4\ln N_a(t-a+a^\star)}{N_a(t-a+a^\star)}}\,.
\]
All in all, we proved the inclusion: for $t \geq K+1$,
\begin{align*}
\{ A_t = a \}
\subseteq & \qquad
\left\{ A_t = a \ \ \mbox{and} \ \ \widehat{\mu}_{a,t-1} - \mu^\star >
- \sqrt{\frac{4\ln N_a(t-1)}{N_a(t-1)}} \, \right\} \\
& \cup
\left\{ A_{t-a+a^\star} = a^\star \ \ \mbox{and} \ \ \widehat{\mu}_{a^\star,t-a+a^\star} - \mu^\star \leq
-\sqrt{\frac{4\ln N_a(t-a+a^\star)}{N_a(t-a+a^\star)}} \, \right\}.
\end{align*}
We now only sketch the next argument, as we proceed similarly to all multi-armed bandit analyses, by resorting to Doob's optional sampling theorem,
which asserts that the rewards $Y_s$ obtained at those rounds $s$ when $A_s = a$ are independent
and identically distributed according to $\nu_a$. We denote by $\overline{\mu}_{a,n}$ the
empirical average of the first $n$ rewards obtained by arm $a$ during the game. Then,
\begin{align}
\nonumber
\E_{\unu} \! \bigl[ N_a(T) \bigr] & \leq 1 + \sum_{t=K+1}^T \P\!\left\{ A_t = a \ \ \mbox{and} \ \ \widehat{\mu}_{a,t-1} - \mu^\star >
- \sqrt{\frac{4\ln N_a(t-1)}{N_a(t-1)}} \, \right\} \\
\nonumber
& \qquad + \sum_{t=K+1}^T \P\left\{ A_{t-a+a^\star} = a^\star \ \ \mbox{and} \ \ \widehat{\mu}_{a^\star,t-a+a^\star} - \mu^\star \leq
-\sqrt{\frac{4\ln N_a(t-a+a^\star)}{N_a(t-a+a^\star)}} \, \right\} \\
\label{eq:ub-end1}
& \leq 1 + \sum_{n \geq 1} \, \P\!\left\{ \overline{\mu}_{a,n} - \mu^\star >
- \sqrt{\frac{4\ln n}{n}} \, \right\} + \sum_{n \geq 1} \, \P\!\left\{ \overline{\mu}_{a^\star,n} - \mu^\star \leq
- \sqrt{\frac{4\ln n}{n}} \, \right\}.
\end{align}
As indicated already in~\citemor{BPR},
for each arm $a$,
the sub-Gaussian assumption on $\nu_a$, together with a Cr{\'a}mer--Chernoff bound, indicates that
for all $n \geq 1$ and all $\varepsilon > 0$,
\begin{equation}
\label{eq:Cramer}
\max\Bigl\{
\P\bigl\{ \overline{\mu}_{a,n} - \mu_a \geq \varepsilon \}, \,\,
\P\bigl\{ \overline{\mu}_{a,n} - \mu_a \leq -\varepsilon \} \Bigr\}
\leq \exp \biggl( - \frac{n \varepsilon^2 }{ 2 } \biggr)\,.
\end{equation}
We substitute this inequality in the bound~\eqref{eq:ub-end1} obtained above.
On the one hand, for $a^\star$,
\begin{equation}
\sum_{n \geq 1} \, \P\!\left\{ \overline{\mu}_{a^\star,n} - \mu^\star \leq
- \sqrt{\frac{4\ln n}{n}} \, \right\}
\leq \sum_{n \geq 1} \, n^{-2} \leq 2\,.
\label{eq:first_term_optimal_arm}
\end{equation}
On the other hand, for $a$, we rewrite $\mu^\star = \mu_a + \Delta_a$ and get
\begin{align*}
\sum_{n \geq 1} \, \P\!\left\{ \overline{\mu}_{a,n} - \mu^\star >
- \sqrt{\frac{4\ln n}{n}} \, \right\}
& = \sum_{n \geq 1} \, \P\!\left\{ \overline{\mu}_{a,n} - \mu_a >
\Delta_a - \sqrt{\frac{4\ln n}{n}} \, \right\}.
\end{align*}
To upper bound the latter sum, we denote by $n_0$ the smallest integer $k \geq 3$, if it exists, such that:
\begin{equation}
\label{eq:n0}
\Delta_a-\sqrt{\frac{4\ln k}{k}}\geq\frac{\Delta_a}{2}\,, \qquad \mbox{that is}, \qquad
\sqrt{\frac{4\ln k}{k}} \leq \frac{\Delta_a}{2}\,.
\end{equation}
As $x \mapsto \sqrt{(\ln x)/x}$ is decreasing on $[3,+\infty)$, we have
\[
\forall n \geq n_0, \qquad
\Delta_a-\sqrt{\frac{4\ln n}{n}}\geq\frac{\Delta_a}{2}\,,
\]
and thus
\[
\sum_{n \geq 1} \, \P\!\left\{ \overline{\mu}_{a,n} - \mu_a >
\Delta_a - \sqrt{\frac{4\ln n}{n}} \, \right\}
\leq n_0-1 + \sum_{n \geq n_0} \, \P\!\left\{ \overline{\mu}_{a,n} - \mu_a >
\frac{\Delta_a}{2} \, \right\}\,.
\]
Note that the above inequality also holds with $n_0 = 2$ when no $k \geq 3$ satisfies~\eqref{eq:n0}.
We use \eqref{eq:Cramer} and
a comparison to an integral to get
\[
\sum_{n \geq n_0} \, \P\!\left\{ \overline{\mu}_{a,n} - \mu_a >
\frac{\Delta_a}{2} \, \right\}
\leq \sum_{n \geq n_0} \, \exp\biggl( - \frac{ n \Delta_a^2 }{ 8 }\biggr)
\leq \bigintsss_{n_0-1}^{+\infty} \exp\biggl( - \frac{ x \Delta_a^2} { 8 }\biggr) \,\mathrm{d} x
\leq \frac{8}{\Delta_a^2}\,.
\]
Substituting the above bounds and~\eqref{eq:first_term_optimal_arm} into~\eqref{eq:ub-end1},
we showed so far that
\[
\E_{\unu} \! \bigl[ N_a(T) \bigr] \leq n_0 + 2 + \frac{8}{\Delta_a^2}\,.
\]
The proof is concluded by upper bounding $n_0$, based on~\eqref{eq:n0}.
If $\Delta_a \leq 4 \sqrt{(\ln 3)/3}$, then the $n_0$ defined in~\eqref{eq:n0} exists.
In this case, we denote by $x_0 \in [3,+\infty)$ the real number such that
\[
\sqrt{\frac{4\ln x_0}{x_0}} \leq \frac{\Delta_a}{2}
\qquad \mbox{that is}, \qquad
x_0 = \frac{16 \ln x_0}{\Delta_a^2}\,.
\]
We have $n_0 = \lceil x_0 \rceil \leq x_0 + 1$.
Since
\[
x_0 = \frac{16 \ln x_0}{\Delta_a^2} = \frac{32 \ln\bigl(4/\Delta)}{\Delta_a^2}
+ \frac{16}{\Delta_a^2} \ln \bigl( \ln x_0 \bigr)\,,
\]
we suspect that $x_0$ should not be too much larger than $32 \ln\bigl(4/\Delta) \big/ \Delta_a^2$.
Indeed, using the inequality $\ln(u) \leq u$, we see that
\[
x_0 = \frac{16 \ln x_0}{\Delta_a^2}
= \frac{160 \ln x_0^{1/10}}{\Delta_a^2}
\leq \frac{160 \, x_0^{1/10}}{\Delta_a^2}\,,
\qquad \mbox{thus} \qquad
x_0 \leq \left( \frac{160}{\Delta_a^2} \right)^{\!\! 10/9}\,.
\]
Therefore,
\[
x_0 = \frac{16 \ln x_0}{\Delta_a^2}
\leq \frac{16}{\Delta_a^2} \ln \! \left( \frac{160}{\Delta_a^2} \right)^{\!\! 10/9}
\leq \frac{16 \times (10/9) \times 2}{\Delta_a^2} \ln \frac{13}{\Delta^2}
\leq \frac{36}{\Delta_a^2} \ln \frac{13}{\Delta^2}\,.
\]
When the $n_0$ defined in~\eqref{eq:n0} does not exist and we take $n_0 = 2$,
we may still bound $n_0$ by $1$ plus the bound above on $x_0$ (as the latter is larger than~$1$).
The theorem follows, after substitution of all the bounds, together with the inequality
$8 \leq 36 \ln(17) - 36 \ln(13)$.
\end{MORproof}

\end{APPENDICES}

\section*{Acknowledgments.}

The authors thank S\'ebastien Gerchinovitz and Vianney Perchet
for stimulating discussions and comments. They are grateful to the anonymous associate
editor and reviewers for their thoughtful feedback and remarks.

This work was partially supported by the CIMI (Centre International de Math{\'e}matiques et d'Informatique)
Excellence program while Gilles Stoltz visited Toulouse in November 2015.
The authors acknowledge the support of the French Agence Nationale de la Recherche (ANR), under grants
ANR-13-BS01-0005 (project SPADRO) and ANR-13-CORD-0020 (project ALICIA).
Gilles Stoltz would like to thank Investissements d'Avenir (ANR-11-IDEX-0003/Labex Ecodec/ANR-11-LABX-0047) for financial support.

\bibliographystyle{ormsv080}
\bibliography{biblio-BLB}

\end{document}